\documentclass[review,hidelinks,onefignum,onetabnum]{siamart250211}

\nolinenumbers
\usepackage{float}
\usepackage{subcaption}
\usepackage{multirow} 
\usepackage{array}      
\usepackage{booktabs}
\usepackage[skip=2pt]{caption}
\usepackage{xr-hyper} 
\usepackage{hyperref}
\usepackage{caption}
\usepackage[utf8]{inputenc}
\usepackage{bm}

\usepackage{lipsum}
\usepackage{amsfonts}
\usepackage{graphicx}
\usepackage{epstopdf}
\usepackage{algorithmic}
\ifpdf
  \DeclareGraphicsExtensions{.eps,.pdf,.png,.jpg}
\else
  \DeclareGraphicsExtensions{.eps}
\fi

\newsiamremark{remark}{Remark}
\newsiamremark{assumption}{Assumption}
\newsiamremark{hypothesis}{Hypothesis}
\crefname{hypothesis}{Hypothesis}{Hypotheses}
\newsiamthm{claim}{Claim}
\newsiamremark{fact}{Fact}
\crefname{fact}{Fact}{Facts}

\headers{}{}

\title{Two-Sided Projection-Based Model Reduction of Parametric Generalized Algebraic Riccati Equations via Grassmann-Manifold Optimization\thanks{Submitted to the editors DATE.
\funding{This work was supported in part by the National Natural Science Foundation of China (12471405, 12671450), the National Key Research and Development Program of China (2023YFB3001604),
and the Science and Technology Innovation Program of Hunan Province (2025RC3080).}}}

\author{Xuelian Wen\thanks{Key Laboratory of Intelligent Computing and Information Processing of Ministry of Education, Key Laboratory for Computation and Simulation in Science and Engineering, the School of Mathematics and Computational Science, Xiangtan University, Xiangtan, Hunan, China
  (\email{202421511236@smail.xtu.edu.cn}, \email{zhangjuan@xtu.edu.cn})}.
\and Qiuqi Li\thanks{the School of Mathematics, Hunan University, Changsha, Hunan, China
  (\email{qli28@hnu.edu.cn})}.
\and Juan Zhang\footnotemark[2]}

\usepackage{amsopn}

\ifpdf
\hypersetup{
	pdftitle={Two-Sided Projection-Based Model Reduction of Parametric  Riccati Equations via Grassmann-Manifold Optimization},
	pdfauthor={Xuelian Wen, Qiuqi Li, Juan Zhang}
}
\fi

\usepackage{xr}

\newcommand{\R}{\mathbb{R}}
\newcommand{\norm}[1]{\left\| #1 \right\|}

\newcommand{\cP}{\mathcal{P}}

\newcommand{\cM}{\mathcal{M}}

\begin{document}
	\maketitle
	
	\begin{abstract}
        Parametric generalized algebraic Riccati equations arise in many-query problems in optimal control, state estimation, and filtering, where repeated full-order solutions may become computationally prohibitive for large-scale systems. In this paper, within a projection-based model reduction framework, we propose a two-sided projection approach that constructs reduced bases by minimizing the aggregate relative two-sided projection error of Riccati solution snapshots. Under suitable dissipativity assumptions, the resulting Galerkin-projected PGARE admits a unique symmetric positive semidefinite stabilizing solution. For error assessment, we employ an existing residual-based a posteriori error estimator for the full-order Riccati residual. A projection-based error decomposition developed in this work provides the theoretical motivation for the proposed basis construction. To accommodate clustered solution manifolds, we extend the framework to localized model reduction by combining parameter-space proximity and solution-space correlation for clustering, with classifier-based online basis selection. Numerical experiments on parameterized heat-conduction problems demonstrate the effectiveness of the proposed global and localized approaches in terms of approximation accuracy and computational efficiency.
	\end{abstract}
	
	\begin{keywords}
		model order reduction, parametric generalized algebraic Riccati equations, Grassmann manifold optimization,  localized reduced basis.
	\end{keywords}
	
	\begin{MSCcodes}
		15A24, 49M30, 65F30, 93C05
	\end{MSCcodes}
	
	\section{Introduction}\label{sec:intro}
	
	\subsection{Background}
	
	Algebraic Riccati equations arise in a wide range of problems in
	systems and control theory, including optimal control, state
	estimation, filtering, and robust control
	\cite{AbouKandil2003,AndersonMoore1990}. When the underlying dynamical system
	depends on physical, geometric, or operating parameters, the
	corresponding Riccati equation becomes parameter dependent and must
	typically be solved repeatedly for different parameter values.
	Such problems arise naturally in design optimization, uncertainty
	quantification, real-time control, and other many-query settings.
	Related parameter-dependent matrix equations also occur in low-gain
	feedback design \cite{Zhou2008Parametric}, parametric model reduction
	\cite{Son2021}, and consensus problems for multi-agent systems
	\cite{Kim2009}. 
	
	In this paper, we consider the parameter-dependent generalized
	algebraic Riccati equation (PGARE)
	\begin{equation}
		\mathcal{R}(X;\mu)
		:=
		A(\mu)^\top X E(\mu)
		+
		E(\mu)^\top X A(\mu)
		-
		E(\mu)^\top X G(\mu) X E(\mu)
		+
		Q(\mu)
		=
		0,
		\label{eq:PGARE}
	\end{equation}
	where
	\(\mu\in\mathcal{P}\subset\mathbb{R}^d\) is a parameter vector,
	\[
	G(\mu)
	=
	B(\mu)B(\mu)^\top
	\succeq 0,
	\qquad
	Q(\mu)
	=
	C(\mu)^\top C(\mu)
	\succeq 0.
	\]
	We assume that the coefficient matrices depend affinely on the
	parameter, namely,
	\begin{equation}
		\begin{aligned}
			A(\mu)
			&=
			\sum_{i=1}^{n_A}
			\theta_i^A(\mu)A_i,
			&
			E(\mu)
			&=
			\sum_{i=1}^{n_E}
			\theta_i^E(\mu)E_i,
			\\
			B(\mu)
			&=
			\sum_{i=1}^{n_B}
			\theta_i^B(\mu)B_i,
			&
			C(\mu)
			&=
			\sum_{i=1}^{n_C}
			\theta_i^C(\mu)C_i,
		\end{aligned}
		\label{eq:affine_decomposition}
	\end{equation}
	where
	\[
	A_i,E_i\in\mathbb{R}^{N\times N},
	\qquad
	B_i\in\mathbb{R}^{N\times n},
	\qquad
	C_i\in\mathbb{R}^{p\times N},
	\]
	with \(n,p\ll N\).
		
	In practical applications, particularly in optimal control, state estimation, and filtering, one is typically interested not in an arbitrary solution of \eqref{eq:PGARE}, but rather in the unique \emph{stabilizing solution}, denoted by $X^*(\mu)$. 
	This specific solution is characterized by two key properties: it is symmetric positive semidefinite, i.e., $X^*(\mu) = X^*(\mu)^\top \succeq 0$, and it guarantees the stability of the corresponding closed-loop matrix pencil $(E(\mu), A(\mu) - G(\mu)X^*(\mu)E(\mu))$. 
	Under standard stabilizability and detectability assumptions, the PGARE admits a unique stabilizing solution for all $\mu\in \mathcal{P}$.

	\subsection{Existing approaches and motivation}
	
	A substantial literature is devoted to the numerical solution of
	algebraic Riccati equations. For moderate problem sizes, classical
	direct approaches include Schur-based methods
	\cite{Laub1979} and structure-preserving algorithms
	\cite{ChuFanLin2005}. For large-scale problems, iterative and
	low-rank techniques have been developed, including Newton-type methods
	\cite{BennerByers1998,Feitzinger2009}, low-rank ADI methods
	\cite{Benner2016,BennerLiPenzl2008}, rational Krylov subspace methods
	\cite{Druskin2011}, and matrix-sign-function approaches
	\cite{Byers1987}; see, for example,
	\cite{BennerSaak2013} for surveys and further
	references. Projection-based system-reduction techniques, such as
	balanced truncation
	\cite{GugercinAntoulas2004,Stykel2004,Moore1981}
	and interpolatory \(\mathcal{H}_2\)-optimal reduction
	\cite{GugercinAntoulasBeattie2008,AntoulasBeattieGugercin2020},
	provide additional tools for reducing large dynamical systems.
	
	For parameter-dependent problems, parametric model order reduction
	(PMOR) provides an offline-online strategy for avoiding repeated
	large-scale computations
	\cite{Benner2015Survey,BennerGrivetTalalal2005}.
	Among these approaches, reduced basis methods have been extensively
	developed for parametrized problems. Originally introduced and
	systematically analyzed for parametrized partial differential
	equations
	\cite{HesthavenRozzaStamm2016,RozzaHuynhPatera2008},
	the reduced basis methodology has also been extended to
	parameter-dependent matrix equations. Son and Stykel
	\cite{Son2017} developed reduced basis methods for
	parameter-dependent Lyapunov equations together with rigorous
	a posteriori error estimates. Schmidt and Haasdonk
	\cite{SchmidtHaasdonk2018} considered large-scale
	parameter-dependent algebraic Riccati equations and constructed
	reduced spaces using a residual-driven greedy strategy.
	
	The present work adopts a different basis-construction principle.
	Let
	\[
	\left\{
	X^*(\mu_j)
	\right\}_{j=1}^{s}
	\]
	denote a collection of symmetric Riccati solution snapshots associated
	with training parameters
	\(\{\mu_j\}_{j=1}^{s}\subset\mathcal{P}\).
	For an orthonormal basis
	\(V\in\mathbb{R}^{N\times r}\),
	the reduced approximation is represented in the two-sided form
	\[
	\widetilde X(\mu)
	=
	V\widehat X(\mu)V^\top.
	\]
	This representation suggests that the quality of the reduced space should be assessed through a two-sided projection rather than through a conventional one-sided snapshot approximation. This two-sided projection criterion is exactly the objective underlying the two-dimensional singular value decomposition (2D-SVD) introduced by Ding and Ye~\cite{DingY05} for 2D maps and images, where the reconstruction quality of a data matrix is measured after projecting it from both sides, rather than by the one-sided action $VV^\top$ alone. Since the stabilizing solution $X^*(\mu)$ is symmetric positive semidefinite for all $\mu \in \mathcal{P}$, projecting from both sides with the same orthonormal basis $V$ preserves this symmetry; that is, $V^\top X^* V = (V^\top X^* V)^\top \succeq 0$. There is therefore no need to optimize separate left and right bases $V_L, V_R$, and the 2D-SVD objective reduces to the single-subspace problem considered here. We therefore seek $V$ by minimizing
	\[
	\min_{V^\top V=I_r}
	\sum_{j=1}^{s}
	\frac{
		\bigl\|
		X^*(\mu_j)
		-
		VV^\top X^*(\mu_j)VV^\top
		\bigr\|_F^2
	}{
		\|X^*(\mu_j)\|_F^2
	}.
	\]
	In contrast to greedy reduced basis constructions, which enlarge the
	reduced space sequentially according to an error or residual
	criterion, the above formulation optimizes a prescribed
	\(r\)-dimensional subspace over all training snapshots simultaneously.
	
	A Pythagorean identity for the two-sided orthogonal projection shows
	that the preceding minimization is equivalent to maximizing
	\[
	\sum_{j=1}^{s}
	\frac{
		\|V^\top X^*(\mu_j)V\|_F^2
	}{
		\|X^*(\mu_j)\|_F^2
	},
	\]
	which measures the aggregate Frobenius energy captured by the reduced
	space. The orthogonality constraint naturally leads to an optimization
	problem on the Grassmann manifold
	\(\mathrm{Gr}(r,N)\). Riemannian optimization methods on matrix
	manifolds are well established
	\cite{AbsilMahonySepulchre2008,Boumal2023,manopt};
	here they provide an efficient computational framework for optimizing
	the Riccati-specific two-sided projection criterion.

	\subsection{Contributions}
	
	The main contributions of this paper are summarized as follows.
	
	\begin{enumerate}
		
		\item
		\textbf{Riccati-oriented reduced-space construction by Grassmann optimization of the two-sided projection.}
		We formulate the reduced-space construction for the PGARE as the
		minimization of an aggregate relative two-sided projection error of
		Riccati solution snapshots. An exact Pythagorean identity transforms
		this problem into the maximization of the captured Frobenius energy
		subject to an orthogonality constraint. The resulting nonconvex
		optimization problem is solved on the Grassmann manifold by a
		Riemannian conjugate-gradient method. 
		
		\item
		\textbf{Theoretical motivation for the optimized basis via projection-based error decomposition.}
		To guarantee the accuracy of reduced approximation, we recall a residual-based error estimator. 
		We then provide a theoretical motivation for the manifold-optimized basis 
		construction through a projection-based error decomposition. Rather than 
		constructing an independent error estimator, we show how the projection 
		component, controlled through the Lipschitz continuity of the 
		parameter-to-solution map, relates to the approximation quality of 
		the reduced space. This bridges the gap between the proposed 
		subspace optimization and the residual-based verification framework.
		
		\item
		\textbf{Localized reduced spaces for clustered solution manifolds.}
		To improve the efficiency of a single global reduced space when the
		Riccati solution manifold exhibits a clustered geometry, we extend the
		framework to a localized setting. A hybrid dissimilarity measure
		combining parameter-space proximity and the Frobenius-angle
		dissimilarity of the solution snapshots is used to partition the
		training set by \(k\)-medoids clustering
		\cite{KaufmanRousseeuw1990}. A separate manifold-optimized reduced
		basis is then constructed for each cluster, and an online classifier
		selects the corresponding local reduced model for a new parameter.
		The residual-based error estimator can be applied on a clusterwise basis,
		while the active online reduced dimension can remain substantially
		smaller than that required by a single global basis.
		
	\end{enumerate}

	\subsection{Organization}
	The remainder of the paper is organized as follows.
	Section~\ref{sec:framework} develops the manifold basis construction and projection-based reduction
	framework.
	 Section~\ref{sec:error} recalls the residual-based error estimator for algebraic Riccati equations from the literature and provides the theoretical motivation for the optimized basis via projection-based error decomposition. Section~\ref{sec:numerics} reports numerical experiments. Section~\ref{sec:conclusion} draws conclusions and
	outlines directions for future work.

	\section{Projection-Based PGARE Reduction via Grassmann-Manifold Optimization}
	\label{sec:framework}
	
	In this section, we adopt a projection-based model reduction framework for the PGARE. The key idea is to first construct a basis by solving an optimization problem on the Grassmann manifold. The original equation is then projected onto this subspace to obtain a reduced-order equation that can be solved efficiently for each parameter. We first describe the projection framework and the reduced-order equation, and then present approaches for constructing the basis: a global basis for the entire parameter domain and a localized variant for clustered solution manifolds.
	
	\subsection{Projection-Based Reduction of PGARE}
	\label{subsec:rom}
	
	Let
	\[
	V\in\mathbb{R}^{N\times r},
	\
	V^\top V=I_r,
	\
	r\ll N
	\]
	be a reduced basis matrix. We seek an approximation to the
	full-order stabilizing solution \(X^*(\mu)\) in the form
	\[
	\widetilde X(\mu)
	=
	V\widehat X(\mu)V^\top,
	\]
	where
	\(\widehat X(\mu)\in\mathbb{R}^{r\times r}\).
	Substituting this approximation into the full-order residual and
	imposing the Galerkin condition
	\[
	V^\top
	\mathcal{R}
	\bigl(
	V\widehat X(\mu)V^\top;\mu
	\bigr)
	V
	=
	0
	\]
	gives the reduced PGARE
	\begin{equation}
		\widehat{\mathcal{R}}(\widehat X;\mu)
		:=
		\widehat A(\mu)^\top
		\widehat X
		\widehat E(\mu)
		+
		\widehat E(\mu)^\top
		\widehat X
		\widehat A(\mu)
		-
		\widehat E(\mu)^\top
		\widehat X
		\widehat G(\mu)
		\widehat X
		\widehat E(\mu)
		+
		\widehat Q(\mu)
		=
		0.
		\label{eq:reduced_PGARE}
	\end{equation}
	The reduced coefficient matrices are defined by
	\[
	\begin{aligned}
		\hat{A}(\mu) &= V^\top A(\mu) V = \sum_{i=1}^{n_A} \theta_i^A(\mu) V^\top A_i V = \sum_{i=1}^{n_A} \theta_i^A(\mu) \hat{A}_i,\\
		\hat{E}(\mu) &= V^\top E(\mu) V = \sum_{i=1}^{n_E} \theta_i^E(\mu) V^\top E_i V = \sum_{i=1}^{n_E} \theta_i^E(\mu) \hat{E}_i,\\
		\hat{B}(\mu) &= V^\top B(\mu) = \sum_{i=1}^{n_B} \theta_i^B(\mu) V^\top B_i = \sum_{i=1}^{n_B} \theta_i^B(\mu) \hat{B}_i,\\
		\hat{C}(\mu) &= C(\mu) V = \sum_{i=1}^{n_C} \theta_i^C(\mu) C_i V = \sum_{i=1}^{n_C} \theta_i^C(\mu) \hat{C}_i,
	\end{aligned}
	\]
	and consequently
	\[
	\hat{G}(\mu) = \hat{B}(\mu) \hat{B}(\mu)^\top, \qquad \hat{Q}(\mu) = \hat{C}(\mu)^\top \hat{C}(\mu).
	\]
	The reduced equation \eqref{eq:reduced_PGARE} is of the same form as the original PGARE \eqref{eq:PGARE}, but of a much smaller dimension $r \ll N$, so it can be solved efficiently for each parameter $\mu$ using standard solvers for algebraic Riccati equations.
	
	The following lemma provides a sufficient structural condition that ensures the reduced PGARE is well-posed and admits a unique symmetric positive semidefinite stabilizing solution.
	\begin{lemma}[Inheritance of dissipativity under Galerkin projection]
		\label{lem:dissipative}
		Assume that for all $\mu \in \mathcal P$,
		\begin{equation}
			E(\mu) = E(\mu)^\top \succ 0, \ A(\mu) + A(\mu)^\top \prec 0. \nonumber
		\end{equation}
		Then for any reduced basis $V\in\mathbb{R}^{N\times r}$ with $V^\top V = I_r$, the projected matrices
		\[
		\hat E(\mu) = V^\top E(\mu) V, \ \hat A(\mu) = V^\top A(\mu) V
		\]
		satisfy
		\[
		\hat E(\mu) \succ 0, \ \hat A(\mu) + \hat A(\mu)^\top \prec 0.
		\]
		Consequently, the reduced pencil $\lambda \hat E(\mu) - \hat A(\mu)$ is stable, and the reduced system $(\hat E(\mu), \hat A(\mu), \hat B(\mu))$ is stabilizable and $(\hat E(\mu), \hat A(\mu), \hat C(\mu))$ is detectable. By the classical Riccati theory for descriptor systems, the reduced PGARE \eqref{eq:reduced_PGARE} admits a unique symmetric positive semidefinite stabilizing solution $\hat X(\mu)$.
		
		\begin{proof}
			The positive definiteness of $\hat E(\mu)$ follows directly from $E \succ 0$ and $V$ having full column rank. For the dissipativity of $\hat A(\mu)$, take any nonzero $y \in \mathbb R^r$ and set $x = Vy \neq 0$. Then
			\[
			y^\top(\hat A(\mu) + \hat A^\top(\mu))y = y^\top V^\top (A(\mu) + A^\top(\mu)) V y = x^\top (A(\mu) + A^\top(\mu)) x < 0.
			\]
			Thus $\hat A(\mu) + \hat A^\top(\mu) \prec 0$. Stability of the reduced pencil follows from standard generalized eigenvalue theory. Since the reduced system is open-loop stable, stabilizability and detectability are automatic. The existence of a unique SPSD stabilizing solution follows from the classical Riccati theory for descriptor systems \cite{LancasterRodman1995}.
		\end{proof}
	\end{lemma}
	
	The reduced-order solution described above is independent of how the projection matrix $V$ is constructed. In the following, we present two strategies for obtaining $V$.

	\subsection{A global subspace construction for PGARE on a Grassmann manifold}
	\label{sec:global}
	
	Let
	\begin{equation}
		\mathcal{P}_{\mathrm{train}}=\{\mu_1,\ldots,\mu_{n_s}\}\subset\mathcal{P}
		\label{eq:training_set}
	\end{equation}
	be a finite training set, and let
	\begin{equation}
		X_j:=X^*(\mu_j),\qquad j=1,\ldots,n_s,
		\label{eq:snapshots}
	\end{equation}
	be the corresponding stabilizing solution snapshots.  We seek a single
	global subspace of dimension $r$ that represents the snapshots accurately
	under two-sided projection.

	For $V^\top V=I_r$, define
	\begin{equation}
		\Phi(V)
		:=\sum_{j=1}^{n_s}
		\frac{\norm{X_j-VV^\top X_jVV^\top}_F^2}{\norm{X_j}_F^2}.
		\label{eq:projection_objective}
	\end{equation}
	The normalization prevents snapshots of large Frobenius norm from dominating
	the basis construction. While this objective is evaluated only on the finite training set, Section~\ref{subsec:basis_motivation} establishes a direct link between $\Phi(V)$ and the projection error indicator $f(\mu)$, showing that optimizing the Grassmann objective also improves the subspace approximation quality at
    unseen parameters.
	
	\begin{proposition}[Pythagorean identity]
		\label{prop:pythagorean}
		For any $X\in\mathbb{R}^{N\times N}$ and any orthogonal projector $P=VV^\top$,
		\begin{equation}
			\norm{X-PXP}_F^2
			=\norm{X}_F^2-\norm{V^\top XV}_F^2.
			\label{eq:pythagorean_identity}
		\end{equation}
		Consequently, minimizing \eqref{eq:projection_objective} is equivalent to
		maximizing
		\begin{equation}
			J(V):=\sum_{j=1}^{n_s}w_j\norm{V^\top X_jV}_F^2,
			\qquad
			w_j:=1/\norm{X_j}_F^{2}.
			\label{eq:captured_energy}
		\end{equation}
	\end{proposition}
	
	\begin{proof}
		The Frobenius inner product satisfies
		\[
		\langle X, PXP\rangle
		=\operatorname{tr}(X^\top P XP)
		=\operatorname{tr}((V^\top XV)^\top(V^\top XV))
		=\norm{V^\top XV}_F^2.
		\]
		Moreover, $\norm{P XP}_F=\norm{V^\top XV}_F$.  Expanding the square of
		$\norm{X-P XP}_F$ proves \eqref{eq:pythagorean_identity}.  Summing the
		normalized identities gives the equivalence.
	\end{proof}
	
\eqref{eq:captured_energy} directly measures
	the energy captured by the two-sided Riccati reconstruction.  It also reveals
	the correct optimization geometry.
	
	\begin{proposition}[Grassmann invariance]
		\label{prop:grassmann_invariance}
		For every orthogonal $O\in\mathbb{R}^{r\times r}$,
		\begin{equation}
			J(VO)=J(V).
			\label{eq:grassmann_invariance}
		\end{equation}
		Hence $J$ depends only on $\operatorname{range}(V)$ and defines a smooth function on the
		Grassmann manifold $\operatorname{Gr}(r,N)$.
	\end{proposition}
	
	\begin{proof}
		For every snapshot,
		\[
		\norm{(VO)^\top X_j(VO)}_F
		=\norm{O^\top(V^\top X_jV)O}_F
		=\norm{V^\top X_jV}_F.
		\]
		Summation proves the claim.
	\end{proof}
	
	
	Because every $X_j$ is symmetric, differentiation of
	\eqref{eq:captured_energy} gives the Euclidean gradient
	\begin{equation}
		\nabla J(V)
		=4\sum_{j=1}^{n_s}w_jX_jV(V^\top X_jV).
		\label{eq:euclidean_gradient}
	\end{equation}
	For the canonical embedded representation of $\operatorname{Gr}(r,N)$, its horizontal
	projection is
	\begin{equation}
		\operatorname{grad}J(V)
		=(I_N-VV^\top)\nabla J(V).
		\label{eq:riemannian_gradient}
	\end{equation}
	Optimization problem \eqref{eq:captured_energy} allows
	the use of standard Riemannian trust-region or conjugate-gradient methods;
	see \cite{AbsilMahonySepulchre2008,Boumal2023,manopt}.
	
	For a large-scale PGARE, the stabilizing snapshots are typically stored in
	low-rank form
	\begin{equation}
		X_j\approx Z_jZ_j^\top,
		\ Z_j\in\mathbb{R}^{N\times k_j},\ k_j\ll N.
		\label{eq:low_rank_snapshots}
	\end{equation}
	Then
	\begin{equation}
		X_jV=Z_j(Z_j^\top V),
		\
		V^\top X_jV=(V^\top Z_j)(V^\top Z_j)^\top,
		\label{eq:low_rank_evaluation}
	\end{equation}
	so neither $X_j$ nor $VV^\top X_jVV^\top$ needs to be formed as a dense
	$N\times N$ matrix.  Moreover,
	\begin{equation}
		\norm{X_j}_F^2
		=\norm{Z_j^\top Z_j}_F^2,
		\label{eq:low_rank_snapshot_norm}
	\end{equation}
	which also makes the weights inexpensive to compute.

	\begin{algorithm}[t]
		\caption{Global Grassmann basis construction}
		\label{alg:global_basis}
		\begin{algorithmic}[1]
			\REQUIRE Training set $\mathcal{P}_{\mathrm{train}}$, target reduced dimension $r$.
			\ENSURE an orthonormal representative $V$ of the optimized subspace.
			\FOR{$j=1,\ldots,n_s$}
			\STATE Compute a low-rank stabilizing snapshot $X_j = Z_j Z_j^\top$.
			\STATE Set $w_j = 1/\norm{Z_j^\top Z_j}_F^{2}$.
			\ENDFOR
			\STATE Compute an orthonormal initial basis $V_0$.
			\STATE Maximize $J$ in \eqref{eq:captured_energy} on the Grassmann manifold $\operatorname{Gr}(r,N)$, using the Riemannian gradient \eqref{eq:riemannian_gradient} and the low-rank evaluation formulas \eqref{eq:low_rank_evaluation}.
			\RETURN an orthonormal representative $V$ of the optimized subspace.
		\end{algorithmic}
	\end{algorithm}

	\subsection{Local subspace construction for PGARE on the Grassmann manifold}
	\label{sec:localized}
	
	The global basis approach described in the previous subsection minimizes the aggregate projection error over the entire training set. When the solution manifold $\cM := \{X^*(\mu) : \mu \in \cP\}$ is well approximated by a single low-dimensional linear subspace, this construction is near-optimal. However, when $\cM$ exhibits a clustered or multi-modal geometry, a single global basis tends to compromise between clusters, and the required dimension $r$ for a given tolerance can grow sharply. In this subsection, we present a localized variant of the basis construction: instead of a single global basis, we partition the training set into clusters and construct a dedicated local basis for each cluster via Grassmann-manifold optimization. The projection step and the reduced-order solve then follow Subsection~\ref{subsec:rom} using the appropriate local basis.
	
	To achieve a partition that reflects both the geometry of the solution manifold and the proximity in the parameter domain, we define a hybrid dissimilarity measure on the space of solutions that combines the Frobenius-angle
	dissimilarity with the parameter-space Euclidean distance:
	\begin{equation}\label{eq:corr}
		d(X_i, X_j; \mu_i, \mu_j) := \alpha \cdot \bar d_{\mathrm{param}}(\mu_i, \mu_j) + (1-\alpha) \cdot \bar d_{\mathrm{cc}}(X_i, X_j),
	\end{equation}
	where $\alpha \in [0,1]$ is a weighting parameter. The parameter-space distance is normalized as
	\[
	\bar d_{\mathrm{param}}(\mu_i, \mu_j) := \frac{\|\mu_i - \mu_j\|_2}{\max_{p,q} \|\mu_p - \mu_q\|_2},
	\]
	and the Frobenius-angle dissimilarity is defined as
	\[
	d_{\mathrm{cc}}(X_i, X_j) := 1 - \frac{\langle X_i, X_j\rangle_F}{\|X_i\|_F \|X_j\|_F} = 1 - \cos\theta_{ij},
	\]
	with $\theta_{ij}$ the angle between $X_i$ and $X_j$ in the Frobenius inner product space, and normalized as
	\[
	\bar d_{\mathrm{cc}}(X_i, X_j) := \frac{d_{\mathrm{cc}}(X_i, X_j)}{\max_{p,q} d_{\mathrm{cc}}(X_p, X_q)}.
	\]
	Both components thus lie in $[0,1]$, so the weight $\alpha$ controls their relative contribution in a balanced manner. The correlation term captures directional dissimilarity between solutions—which is relevant for projection-based approximation—while the parameter term enforces locality in the parameter space, ensuring that nearby parameters are assigned to the same cluster for a smooth online classification.
	
	Given the training solutions $\{X^*(\mu_j)\}_{j=1}^s$ and the distance $d$, we partition the training set into $n_c$ clusters using the $k$-medoids algorithm \cite{KaufmanRousseeuw1990}, which operates directly on pairwise distances and is robust to outliers. Let the resulting clusters be $\mathcal{C}_1, \ldots, \mathcal{C}_{n_c}$, with medoids $\mu_{m_1}, \ldots, \mu_{m_{n_c}}$. For each cluster $\mathcal{C}_\ell$, we construct a dedicated local basis $V_\ell \in \mathrm{St}(r_\ell, N)$ by solving the cluster-specific Grassmann manifold optimization:
	\[
	\max_{V_\ell \in \mathrm{St}(r_\ell, N)} \sum_{\mu_j \in \mathcal{C}_\ell} \frac{\|V_\ell^\top X^*(\mu_j) V_\ell\|_F^2}{\|X^*(\mu_j)\|_F^2}.
	\]
	The local reduced matrices are then precomputed for each cluster $\ell$ as $V_\ell^\top A_i V_\ell$, $V_\ell^\top E_i V_\ell$, $V_\ell^\top B_i$, $C_i V_\ell$, respectively.
	
	After the clusters have been determined in the offline stage, we train a classifier $\phi: \mathcal{P} \to \{1,\ldots,n_c\}$ on the labeled training set $\{(\mu_j, \ell_j)\}_{j=1}^s$, where $\ell_j$ is the cluster assignment of $\mu_j$. The classifier can be chosen according to the problem at hand; for instance, linear discriminant analysis (LDA), $k$-nearest neighbors, support vector machines, or decision trees are all viable options.
	
	In the online stage, for a new parameter $\mu$, the trained classifier $\phi$ predicts the appropriate cluster $\ell = \phi(\mu)$, and the online procedure proceeds exactly as described in Subsection~\ref{subsec:rom}, using $V_\ell$ and the corresponding local reduced matrices.
	
	The localized framework is a natural generalization of the global framework: setting $n_c = 1$ recovers the global basis construction of subsection ~\ref{sec:global}. The localized framework is most beneficial when the solution manifold is clustered, i.e., when the global basis requires a large $r$ to achieve a target accuracy. When the solution manifold is unimodal, the overhead of clustering and classification provides no benefit, and the global framework should be preferred.
	\begin{remark}\label{rem:dist_choice}
		The hybrid dissimilarity measure \eqref{eq:corr} combines the Frobenius-angle dissimilarity, which captures the directional structure of the solution manifold relevant for projection-based reduction, with the parameter-space distance, which enforces locality in the parameter domain to ensure a smooth online classification. The weight $\alpha$ controls the trade-off between these two criteria and is set to a fixed value throughout our experiments. The $k$-medoids algorithm is chosen for its robustness to outliers and its compatibility with arbitrary distance functions.
	\end{remark}

	\begin{algorithm}[htbp]
		\caption{Local Grassmann basis construction}
		\label{alg:local}
		\begin{algorithmic}[1]
			\REQUIRE Training set $\{\mu_j\}_{j=1}^s$, number of clusters $n_c$, local dimensions $\{r_\ell\}_{\ell=1}^{n_c}$, weight $\alpha\in[0,1]$, tolerance $\varepsilon$ for Grassmann optimization.
			\ENSURE Local projection matrices $\{V_\ell\}_{\ell=1}^{n_c}$, local reduced components, and classifier $\phi$.
			\STATE Compute normalized parameter-space distance $\bar{d}_{\mathrm{param}}$ and correlation-coefficient distance $\bar{d}_{\mathrm{cc}}$ between all training pairs.
			\STATE Form hybrid distance $d = \alpha\bar{d}_{\mathrm{param}} + (1-\alpha)\bar{d}_{\mathrm{cc}}$.
			\STATE Partition training set into $n_c$ clusters $\{\mathcal{C}_\ell\}_{\ell=1}^{n_c}$ via $k$-medoids using $d$.
			\FOR{$\ell=1,\dots,n_c$}
			\STATE Solve cluster-specific Grassmann manifold optimization:
			\[
			\max_{V_\ell\in\mathrm{St}(r_\ell,N)}
			\sum_{\mu_j\in\mathcal{C}_\ell}
			\frac{\|V_\ell^\top X^*(\mu_j)V_\ell\|_F^2}{\|X^*(\mu_j)\|_F^2}.
			\]
			\STATE Compute local projected affine components:
			$\hat{A}_i^{(\ell)} = V_\ell^\top A_i V_\ell$, $\hat{E}_i^{(\ell)} = V_\ell^\top E_i V_\ell$,
			$\hat{B}_i^{(\ell)} = V_\ell^\top B_i$, $\hat{C}_i^{(\ell)} = C_i V_\ell$.
			\ENDFOR
			\STATE Train classifier $\phi:\mathcal{P}\to\{1,\dots,n_c\}$ on labeled data $\{(\mu_j,\ell_j)\}_{j=1}^s$.
			\RETURN $\{V_\ell\}$, local reduced components, and $\phi$.
		\end{algorithmic}
	\end{algorithm}
	
	Once the local bases and their reduced components are obtained, the online solution for any new parameter $\mu$ proceeds by first classifying $\ell = \phi(\mu)$, then following the reduced-order solve described in Subsection~\ref{subsec:rom} using $V_\ell$ and its associated reduced matrices.

	\section{A Residual-based Error Estimator and Theoretical Motivation for the Optimized Basis}
	\label{sec:error}
	
	In this section, we first recall a classical residual-based error estimator for 
	algebraic Riccati equations~\cite{SchmidtHaasdonk2018}. We then explain the motivation 
	for constructing the basis via Grassmann manifold optimization. Throughout this section, we may drop the parameter $\mu$ from all $\mu$-dependent quantities 
	for brevity, writing, e.g., $X^*$ for $X^*(\mu)$ and $\widehat X$ for 
	$\widehat X(\mu)$.
	
	\subsection{A Residual-based Error Estimator}
\label{subsec:residual_estimation}

We first recall the residual-based error estimator for algebraic Riccati 
equations, which provides a rigorous and online-computable certificate for 
the reduced approximation. We state the result in the following lemma.

\begin{lemma}[Residual-based a posteriori error estimator\cite{SchmidtHaasdonk2018}]
    \label{thm:residual_certificate}
    Let $\widetilde{X} = \widetilde{X}^\top$ be an approximate solution to the 
    generalized Riccati equation~\eqref{eq:PGARE}. Define the closed-loop matrix, 
    the residual norm, and the generalized Lyapunov operator 
    $\mathcal{L}_{\widetilde{X}} : \mathbb{S}^n \to \mathbb{S}^n$ by
    \begin{equation}
        A_c := A - G\widetilde{X}E, \quad
        \varepsilon := \norm{\mathcal{R}(\widetilde{X})}_F, \quad
        \mathcal{L}_{\widetilde{X}}(H) := A_c^\top H E + E^\top H A_c,
    \end{equation}
    and assume $\mathcal{L}_{\widetilde{X}}$ is invertible. Let
    \begin{equation}
        \gamma := \norm{\mathcal{L}_{\widetilde{X}}^{-1}},
        \quad
        L := 2\norm{E}_2^2 \norm{G}_2.
    \end{equation}
    If the validity criterion
    \begin{equation}
        \vartheta := 4\gamma^2 L\varepsilon < 1
        \label{eq:validity_criterion}
    \end{equation}
    holds, then the unique stabilizing solution $X^*$ of \eqref{eq:PGARE} 
    satisfies
    \begin{equation}
        \norm{X^* - \widetilde{X}}_F \le \eta_{\mathrm{res}}
        := \frac{1 - \sqrt{1 - \vartheta}}{2\gamma L}.
        \label{eq:certified_error_bound}
    \end{equation}
\end{lemma}

Having established the theoretical residual-based error estimator, we now clarify how it is evaluated in practice. The computation of the residual norm $\varepsilon$ and the constants $\gamma$ and $L$ must be performed without constructing any full-order $N \times N$ matrices.

\subsection{Efficient Evaluation of Residual-Based Error Estimator}
\label{subsec:offline_online_residual}

The key to an efficient implementation is an offline--online decomposition that exploits the affine parameter dependence of the system matrices and the low-rank structure of the reduced approximation. We explicitly demonstrate that the online complexity of evaluating the residual-based error estimator is independent of the full-order dimension $N$.

\paragraph{The evaluation of the residual norm}
For the residual norm, we rely on the trace-based offline--online decomposition detailed in \cite[Section~5]{SchmidtHaasdonk2018}. By expanding the squared Frobenius norm $\|\mathcal{R}(\widetilde{X};\mu)\|_F^2 = \operatorname{tr}(\mathcal{R}^\top \mathcal{R})$ and exploiting the affine parameter separability together with the cyclic property of the trace, all full-order matrix products are precomputed offline as small $r \times r$ matrices. The online evaluation then only requires assembling these precomputed components and scales as $O(r^3)$, which is completely independent of $N$.

\paragraph{The evaluation of the constant $L$}
The constant $L = 2\|E\|_2^2\|G\|_2$ depends on the parameter $\mu$. Since $E(\mu)$ admits an affine decomposition and $G(\mu)=B(\mu)B(\mu)^\top$ inherits a finite-term decomposition from the affine parameterization of $B(\mu)$, we precompute offline the spectral norms of the basis matrices $\|E_i\|_2$ and $\|G_i\|_2$, together with the maximal values of the parameter coefficients. This yields a parameter-independent upper bound $L_{\max} \ge L(\mu)$ for all $\mu \in \mathcal{P}$. Alternatively, a tighter parameter-dependent bound can be assembled online with $O(n_E + n_G)$ complexity, avoiding any full-order matrix operations.

\paragraph{The evaluation of the stability constant $\gamma$}
The operator norm $\gamma = \|\mathcal{L}_{\widetilde{X}}^{-1}\|$ is not directly accessible in large-scale settings. Several online-computable approaches for evaluating $\gamma$ are provided in \cite{SchmidtHaasdonk2018}. In this work, we adopt two complementary strategies depending on the problem structure. 

\textit{(i) Estimated bound via closed-loop dissipativity.}
If the closed-loop dissipativity margin
\begin{equation}
    \omega := -\lambda_{\max}\!\left(\operatorname{sym}(E^{-1}A_c)\right)
\end{equation}
is positive, then the inverse Lyapunov operator norm satisfies the classical bound
\begin{equation}
    \gamma \le \overline{\gamma} := \frac{\norm{E^{-1}}_2^2}{2\omega},
    \label{eq:gamma_dissipativity}
\end{equation}
where $\operatorname{sym}(M) := \tfrac{1}{2}(M+M^\top)$ is the symmetric part of a matrix $M$. The bound follows from $\|\mathrm{e}^{-E^{-1}A_c t}\|_2 \le \mathrm{e}^{-\omega t}$, which is exactly what $\omega>0$ guarantees, together with the submultiplicativity of the norm.

While \eqref{eq:gamma_dissipativity} holds for an arbitrary pencil $(E,A)$, it does not by itself yield an online bound: for $E\neq I$, the matrix $\operatorname{sym}(E^{-1}A_c)$ bears no useful spectral relation to $\operatorname{sym}(A_c)$, so $\omega$ can only be evaluated by explicitly forming $E^{-1}A_c$, at a cost of $O(N^3)$. The online evaluation derived below therefore requires $E=I$, as is the case for the second numerical example; the first example has $E\neq I$ and is assessed through \eqref{eq:gamma_dissipativity} directly.

For the specific thermal problem considered in this work, the system is described by $E=I$ and an open-loop matrix of the form $A(\mu) = -K(\mu) - \delta I$ with $K(\mu) \succeq 0$, where $\delta$ is the damping parameter denoted by $\gamma_d$ in Section~4. Under these conditions, the closed-loop matrix is $A_c = A - G\widetilde{X}$. Thus, the dissipativity margin satisfies
\begin{equation}
    \omega = -\lambda_{\max}\!\left(\operatorname{sym}(A - G\widetilde{X})\right)
    \ge \delta + \lambda_{\min}\!\left(\operatorname{sym}(G\widetilde{X})\right)
    \ge \delta - \|G\|_2\,\|\widetilde{X}\|_2,
\end{equation}
where we have used the fact that $\operatorname{sym}(A) = -K - \delta I \preceq -\delta I$ and the inequality $\lambda_{\min}(\operatorname{sym}(G\widetilde{X})) \ge -\|\operatorname{sym}(G\widetilde{X})\|_2 \ge -\|G\widetilde{X}\|_2 \ge -\|G\|_2\|\widetilde{X}\|_2$, in which $\|\operatorname{sym}(\cdot)\|_2 \le \|\cdot\|_2$ and the submultiplicativity of the spectral norm are invoked. Since the reduced approximation is symmetric, $\widetilde{X} = V\hat{X}V^\top$, its spectral norm satisfies $\|\widetilde{X}\|_2 = \|V\hat{X}V^\top\|_2 = \|\hat{X}\|_2$ due to the orthonormality of $V$. Consequently, the certified bound for the considered problem becomes
\begin{equation}
    \omega \ge \delta - \|G\|_2\|\hat{X}\|_2.
\end{equation}

Since $\|G\|_2$ is precomputed offline and $\|\hat{X}\|_2$ is evaluated in the reduced space, this certified bound is computed online with $O(r^3)$ complexity, independent of $N$. It is informative only when the right-hand side of the last inequality is positive; otherwise no certification follows, and the general bound \eqref{eq:gamma_dissipativity} has to be used instead.

\textit{(ii) Estimated bound via projected Lyapunov solve.} 
When $E$ is non-identity or the open-loop dissipativity margin is unavailable, $\gamma$ may alternatively be approximated by a projected Lyapunov method. We precompute offline a Lyapunov basis matrix $V_L \in \mathbb{R}^{N \times r_L}$ ($r_L \ll N$) and the reduced matrices $\hat{A}_{c,i} = V_L^\top A_i V_L$, $\hat{E}_L = V_L^\top E V_L$. Online, the reduced closed-loop matrix $\hat{A}_c(\mu) \in \mathbb{R}^{r_L \times r_L}$ is assembled in $O(r_L^2)$ operations, and the Lyapunov equation
\begin{equation}
    \hat{A}_c(\mu)^\top \hat{H} \hat{E}_L + \hat{E}_L^\top \hat{H} \hat{A}_c(\mu) = -I
\end{equation}
is solved for $\hat{H} \in \mathbb{R}^{r_L \times r_L}$ with $O(r_L^3)$ complexity. The constant $\gamma$ is then approximated by $\|\hat{H}\|_2$\cite{SchmidtHaasdonk2018}. 

When a certified upper bound for $\gamma$ is available, the resulting residual estimator provides a rigorous a posteriori error estimator. When $\gamma$ is approximated through the projected Lyapunov solve, the resulting quantity is used as a practical error estimate. Consequently, the entire residual certificate $(\varepsilon, \gamma, L)$ is evaluated online in $O(r^3 + r_L^3)$, which is independent of the full-order dimension $N$.

Having established an online-computable residual-based error estimator and the efficient offline-online evaluation of the constants $\gamma$ and $L$, we are now ready to address the stabilizing property of the reduced-order solution. 
While Lemma~\ref{thm:residual_certificate} guarantees the accuracy of $\widetilde{X}$, it does not automatically ensure that the closed-loop matrix pair $(E, A_c)$ is stable. The following lemma provides a sufficient condition for determining the stability of the approximate solution $\widetilde{X}$.

\begin{lemma}[Stabilizing property of the approximate solution \cite{SchmidtHaasdonk2018}]
    \label{lem:stabilizing_verification}
    Let $\widetilde{X} = V\hat{X}V^\top$ be the reduced-order solution, and let $\eta_{\mathrm{res}}$ be the residual-based error estimator from Lemma~\ref{thm:residual_certificate}. If the condition
    \begin{equation}
        2\gamma \norm{G}_2 \frac{\norm{E}_2}{\norm{E^{-1}}_2} \eta_{\mathrm{res}} \le 1
        \label{eq:stabilizing_criterion}
    \end{equation}
    is satisfied, then $\widetilde{X}$ is a stabilizing approximate solution for the original PGARE.
\end{lemma}
	\subsection{Theoretical motivation for the optimized basis}
	\label{subsec:basis_motivation}
	
	The residual-based estimator in Subsection~\ref{subsec:residual_estimation} provides an online bound for the reduced error under the validity criterion \eqref{eq:validity_criterion}, but it does not explain how the basis \(V\) should be chosen. In this subsection, we answer this 
	question through a projection-based error analysis. The analysis is not 
	an independent error estimator; rather, it decomposes the total error into a 
	projection component and a subspace approximation component. We show that 
	the projection component depends only on the choice of $V$ and serves as 
	a key quantity controlling the error decomposition. This provides the theoretical 
	motivation for the Grassmann manifold optimization in Section~\ref{sec:framework}, 
	which is designed to minimize this component and thus improve the overall 
	accuracy of the reduced approximation.
	
	Define the error matrix $\Delta := X^* - V\hat{X}V^\top$. By inserting the 
	projection $PX^*P$ with $P = VV^\top$, the error admits the decomposition
	\begin{equation}
		\Delta = \underbrace{(X^* - PX^*P)}_{\delta_{\mathrm{proj}}} 
		+ \underbrace{(PX^*P - V\hat{X}V^\top)}_{\delta_{\mathrm{approx}}}.
		\label{eq:decomp}
	\end{equation}
	
	By the triangle inequality and $\|V\hat{X}V^\top\|_F = \|\hat{X}\|_F$, the 
	relative approximate error satisfies
	\begin{equation}
		\frac{\|\Delta\|_F}{\|V\hat{X}V^\top\|_F}
		\le
		\frac{\|\delta_{\mathrm{proj}}\|_F}{\|\hat{X}\|_F}
		+
		\frac{\|\delta_{\mathrm{approx}}\|_F}{\|\hat{X}\|_F}.
		\label{eq:rel_error}
	\end{equation}
	The two terms on the right-hand side are estimated in the next two subsections.
	
	\subsubsection{Projection error}
	\label{par:projection}
	
	We first estimate the projection term. The key ingredient is the Lipschitz 
	continuity of the solution map $\mu \mapsto X^*(\mu)$ and of the absolute 
	projection error $f(\mu) := \|X^*(\mu) - PX^*(\mu)P\|_F$.
	
	We recall the standing assumptions:
	\begin{itemize}
		\item[(A1)] $\mathcal{P} \subset \mathbb{R}^d$ is compact and convex.
		\item[(A2)] $A(\mu), E(\mu), B(\mu), C(\mu)$ admit affine parameterizations 
		with Lipschitz continuous coefficient functions on $\mathcal{P}$.
		\item[(A3)] $E(\mu)$ is nonsingular for all $\mu \in \mathcal{P}$.
		\item[(A4)] The pencil $\lambda E(\mu) - A(\mu)$ is stable for all 
		$\mu \in \mathcal{P}$.
	\end{itemize}
	
	\begin{lemma}[Lipschitz continuity of the inverse]
		\label{lem:inv_lipschitz}
		Let $\mathcal{P} \subset \mathbb{R}^d$ be compact and let $E(\mu)$ be 
		nonsingular and Lipschitz continuous on $\mathcal{P}$. Then 
		$\mu \mapsto E(\mu)^{-1}$ is Lipschitz continuous on $\mathcal{P}$.
	\end{lemma}
	
	\begin{proof}
		For any $\mu, \nu \in \mathcal{P}$, the resolvent identity gives
		\[
		E(\mu)^{-1} - E(\nu)^{-1} 
		= E(\mu)^{-1}\bigl(E(\nu) - E(\mu)\bigr)E(\nu)^{-1}.
		\]
		Taking the Frobenius norm and using submultiplicativity,
		\[
		\|E(\mu)^{-1} - E(\nu)^{-1}\|_F
		\le \|E(\mu)^{-1}\|_2 \, \|E(\nu) - E(\mu)\|_F \, \|E(\nu)^{-1}\|_2.
		\]
		Since $E(\mu)$ is nonsingular on the compact set $\mathcal{P}$,  $M_E := \max_{\mu \in \mathcal{P}} \|E(\mu)^{-1}\|_2$ is finite. 
		Combined with the Lipschitz continuity of $E(\mu)$ with constant $L_E$, 
		we obtain $\|E(\mu)^{-1} - E(\nu)^{-1}\|_F \le M_E^2 L_E \|\mu - \nu\|_2$.
	\end{proof}
	
	\begin{theorem}[Lipschitz continuity of the PGARE solution]
		\label{thm:lipschitz}
		Under (A1)--(A4), there exists $L_X > 0$ such that
		\begin{equation}
			\|X^*(\mu) - X^*(\mu')\|_F \le L_X \|\mu - \mu'\|_2,
			\qquad \forall \mu, \mu' \in \mathcal{P}.
			\label{eq:lipschitz_X}
		\end{equation}
	\end{theorem}
	
	\begin{proof}
		We apply Sun's perturbation theorem \cite{Sun1998} twice. The proof details are provided in the appendix.
	\end{proof}
	
	\begin{lemma}[Lipschitz continuity of $f$]
		\label{lem:f_lipschitz}
		Under the assumptions of Theorem~\ref{thm:lipschitz}, the absolute 
		projection error $f(\mu) := \|X^*(\mu) - PX^*(\mu)P\|_F$ is Lipschitz 
		continuous on $\mathcal{P}$ with constant $L_f \le L_X$.
	\end{lemma}
	
	\begin{proof}
		The reverse triangle inequality and the contraction property 
		 of the orthogonal projection give
		\[
		|f(\mu) - f(\nu)| \le \|X^*(\mu) - X^*(\nu)\|_F \le L_X \|\mu - \nu\|_2.
		\]
	\end{proof}

	\begin{proposition}[Projection error indicator]
		\label{thm:proj_bound}
		For any $\mu \in \mathcal{P}$, let 
		$\mu_{\mathrm{nn}} = \operatorname*{arg\,min}_{j} \|\mu - \mu_j\|_2$ 
		be the nearest training point. Then
		\begin{equation}
			\|\delta_{\mathrm{proj}}(\mu)\|_F = f(\mu) \le d_{\mathrm{comp}}
			:= f(\mu_{\mathrm{nn}}) + L_f \|\mu - \mu_{\mathrm{nn}}\|_2.
			\label{eq:proj_bound}
		\end{equation}
	\end{proposition}
	
	\begin{proof}
		The proof follows directly from Lemma~\ref{lem:f_lipschitz} and the triangle inequality.
	\end{proof}
	
	\begin{remark}[Practical computation of $L_f$]
		\label{rem:Lf_computation}
		The Lipschitz constant $L_f$ is estimated offline from the training data as
		\[
		\hat{L}_f = \max_{i \neq j} \frac{|f(\mu_i) - f(\mu_j)|}{\|\mu_i - \mu_j\|_2}.
		\]
        We emphasize that $\hat{L}_f$ is an empirical estimate obtained from the finite training set; it serves as a practical indicator rather than a certified rigorous upper bound on the true Lipschitz constant $L_f$.
	\end{remark}
	
	\subsubsection{ Subspace approximation error}
	We now estimate the second term of \eqref{eq:rel_error},
	\[
	\|\delta_{\mathrm{approx}}\|_F = \|PX^*P - V\hat{X}V^\top\|_F 
	= \|\hat{X}^* - \hat{X}\|_F,
	\]
	where $\hat{X}^* := V^\top X^* V$. Substituting 
	$X^* = PX^*P + \delta_{\mathrm{proj}}$ into $\mathcal{R}(X^*;\mu) = 0$ and 
	projecting both sides by $V^\top(\cdot)V$ yields the error equation for 
	$\hat{\Delta} := \hat{X}^* - \hat{X}$:
	\begin{equation}
		\hat{A}_c^\top \hat{\Delta}\hat{E} 
		+ \hat{E}^\top \hat{\Delta}\hat{A}_c
		= -\mathcal{E}_\delta + \hat{E}^\top \hat{\Delta}\hat{G}\hat{\Delta}\hat{E},
		\label{eq:hat_lyapunov}
	\end{equation}
	with $\hat{A}_c := \hat{A} - \hat{G}\hat{X}\hat{E}$ and
	\begin{equation}
		\begin{split}
			\mathcal{E}_\delta 
			= {} & V^\top\bigl(A^\top\delta_{\mathrm{proj}}E 
			+ E^\top\delta_{\mathrm{proj}}A 
			- E^\top\delta_{\mathrm{proj}}GPX^*PE \\
			& - E^\top PX^*PG\delta_{\mathrm{proj}}E 
			- E^\top\delta_{\mathrm{proj}}G\delta_{\mathrm{proj}}E\bigr)V.
		\end{split}
		\label{eq:E_delta}
	\end{equation}
	
	Since $PX^*P = V(\hat{X} + \hat{\Delta})V^\top$, we have
	\begin{equation}
		\|PX^*P\|_F = \|\hat{X}^*\|_F \le \|\hat{X}\|_F + \|\hat{\Delta}\|_F.
		\label{eq:PXP_bound}
	\end{equation}
	
	\begin{lemma}
		\label{lem:Edelta_local}
		Let $d := \|\delta_{\mathrm{proj}}\|_F$, 
		$c_1 := 2\|A\|_2\|E\|_2$, $c_2 := \|E\|_2^2\|G\|_2$, and 
		$C_0 := c_1 + 2c_2\|\hat{X}\|_F$. Then
		\begin{equation}
			\|\mathcal{E}_\delta\|_F \le C_0 d + 2c_2 d\|\hat{\Delta}\|_F + c_2 d^2.
			\label{eq:Edelta_local}
		\end{equation}
	\end{lemma}
	
	\begin{proof}
		Taking the Frobenius norm of each term in \eqref{eq:E_delta} and using 
		\eqref{eq:PXP_bound},
		\[
		\|\mathcal{E}_\delta\|_F 
		\le c_1 d + 2c_2 d(\|\hat{X}\|_F + \|\hat{\Delta}\|_F) + c_2 d^2
		= C_0 d + 2c_2 d\|\hat{\Delta}\|_F + c_2 d^2.
		\]
	\end{proof}
	
	Define the reduced Lyapunov operator 
	$\hat{\mathcal{L}}(Z) := \hat{A}_c^\top Z\hat{E} + \hat{E}^\top Z\hat{A}_c$ 
	and the constants
	\begin{equation}
		\hat{\alpha} := \|\hat{\mathcal{L}}^{-1}\|,
		\
		\hat{\beta} := \|\hat{E}\|_2^2 \|\hat{G}\|_2.
	\end{equation}
	
	\begin{lemma}[Quadratic inequality]
		\label{lem:quad_ineq}
		Let $x := \|\hat{\Delta}\|_F$ and $d := \|\delta_{\mathrm{proj}}\|_F$. Then
		\begin{equation}
			\hat{\alpha}\hat{\beta} x^2 + (2\hat{\alpha} c_2 d - 1)x 
			+ \hat{\alpha}(C_0 d + c_2 d^2) \ge 0.
			\label{eq:quad_ineq}
		\end{equation}
	\end{lemma}
	
	\begin{proof}
		Taking norms in \eqref{eq:hat_lyapunov} and using Lemma~\ref{lem:Edelta_local},
		\[
		x \le \hat{\alpha}(C_0 d + 2c_2 d x + c_2 d^2) + \hat{\alpha}\hat{\beta}x^2,
		\]
		which rearranges to \eqref{eq:quad_ineq}.
	\end{proof}
	
	The quadratic has leading coefficient $a = \hat{\alpha}\hat{\beta} > 0$. To 
	extract a finite upper bound for $x$, we require
	\begin{enumerate}
		\item[(C1)] $2\hat{\alpha} c_2 d < 1$;
		\item[(C2)] $(1 - 2\hat{\alpha} c_2 d)^2 - 4\hat{\alpha}^2\hat{\beta}(C_0 d 
		+ c_2 d^2) > 0$.
	\end{enumerate}
	When (C1) and (C2) hold, the two roots $x_1 < x_2$ of the quadratic are 
	positive and given by
	\begin{equation}
		x_{1,2} = \frac{1 - 2\hat{\alpha} c_2 d 
			\pm \sqrt{(1 - 2\hat{\alpha} c_2 d)^2 
				- 4\hat{\alpha}^2\hat{\beta}(C_0 d + c_2 d^2)}}{2\hat{\alpha}\hat{\beta}}.
		\label{eq:roots}
	\end{equation}
	The next result excludes the larger branch $x \ge x_2$.
	
	\begin{proposition}[Small-branch selection]
		\label{prop:small_branch}
		Assume (C1)--(C2) hold, and let $x_1(d) < x_2(d)$ be the two positive 
		roots given in \eqref{eq:roots}. Suppose the residual certificate 
		$\eta_{\mathrm{res}}$ from Lemma~\ref{thm:residual_certificate} satisfies
		\begin{equation}
			\eta_{\mathrm{res}} < x_2(d).
			\label{eq:branch_test}
		\end{equation}
		Then
		\begin{equation}
			\|\hat{\Delta}\|_F \le x_1(d).
			\label{eq:small_branch}
		\end{equation}
	\end{proposition}
	
	\begin{proof}
		By the error decomposition \eqref{eq:decomp} and the orthogonality 
		$V^\top\delta_{\mathrm{proj}}V = 0$, we have
		\[
		\|\hat{\Delta}\|_F = \|\delta_{\mathrm{approx}}\|_F 
		\le \|\Delta\|_F = \|X^* - \widetilde{X}\|_F.
		\]
		Combining with the residual certificate of 
		Lemma~\ref{thm:residual_certificate},
		\[
		\|\hat{\Delta}\|_F \le \eta_{\mathrm{res}}.
		\]
		By \eqref{eq:branch_test}, this implies 
		$\|\hat{\Delta}\|_F < x_2(d)$. The quadratic inequality 
		\eqref{eq:quad_ineq} therefore rules out the large branch 
		$\|\hat{\Delta}\|_F \ge x_2(d)$, so the small branch 
		\eqref{eq:small_branch} holds.
	\end{proof}
	
	\begin{proposition}
		\label{thm:approx_bound}
        Assume (C1)--(C2) hold, and let $x_1(d) < x_2(d)$ be the two positive 
		roots given in \eqref{eq:roots}. Suppose the residual certificate 
		$\eta_{\mathrm{res}}$ from Lemma~\ref{thm:residual_certificate} satisfies
		\begin{equation}
			\eta_{\mathrm{res}} < x_2(d).
		\end{equation}
		Then,
		\begin{equation}
			\|\hat{\Delta}\|_F \le F(d) 
			:= \frac{2\hat{\alpha}(C_0 d + c_2 d^2)}{1 - 2\hat{\alpha} c_2 d}.
			\label{eq:approx_bound}
		\end{equation}
	\end{proposition}
	
	\begin{proof}
		The exact small root is bounded using $1 - \sqrt{1 - t} \le t$ for 
		$t \in [0, 1)$ with $t := 4\hat{\alpha}^2\hat{\beta}(C_0 d + c_2 d^2)/
		(1 - 2\hat{\alpha} c_2 d)^2$.
	\end{proof}
	
	
	\subsubsection{Theoretical Motivation for the Basis Construction}
	We now synthesize the two components analyzed above to clarify the rationale behind the basis construction, rather than to propose an independent error estimator. 
	The following elementary lemma connects the unknown exact projection error \(d\) to the total error \(\|\Delta\|_F\). 
	This connection reveals the intrinsic mechanism by which the projection component dictates the overall approximation quality, thereby providing the theoretical justification for the Grassmann manifold optimization in Section~\ref{sec:framework}, which is designed to minimize this component.
	
	\begin{lemma}
		\label{lem:compression}
		Suppose the reduced closed-loop matrix $\hat{A}_c$ is asymptotically stable, (C1)--(C2) hold, and let $x_1(d) < x_2(d)$ be the two positive roots given in \eqref{eq:roots}. Suppose further that the residual estimator $\eta_{\mathrm{res}}$ from Lemma~\ref{thm:residual_certificate} satisfies $\eta_{\mathrm{res}} < x_2(d)$. Define
        \begin{equation}
            \varepsilon(z) := \frac{z(1 + 2\hat{\alpha} C_0)}{1 - 2\hat{\alpha} c_2 z},
		    \label{eq:epsilon}
        \end{equation}
	    Then the error satisfies $\|\Delta\|_F \le \varepsilon(d)$.
	\end{lemma}
	
	\begin{proof}
		By the triangle inequality and \eqref{eq:decomp},
		\[
		\|\Delta\|_F \le d + \|\hat{\Delta}\|_F 
		\le d + F(d) = \varepsilon(d).
		\]
	\end{proof}

    Recall that \(d_{comp} = f(\mu_{nn}) + \hat{L}_f \|\mu - \mu_{nn}\|_2\), where the training values \(f(\mu_j)\) and the empirical Lipschitz constant \(\hat{L}_f\) are computed offline as in Proposition~\ref{thm:proj_bound} and Remark~\ref{rem:Lf_computation}, respectively.
	Replacing the unknown $d$ by its computable indicator $d_{\mathrm{comp}}$ 
	from Proposition~\ref{thm:proj_bound} yields the final quantity.
	
	\begin{proposition}
		\label{thm:main}
		Let $\mu \in \mathcal{P}$ and 
		$\mu_{\mathrm{nn}} = \operatorname*{arg\,min}_j \|\mu - \mu_j\|_2$. Under 
		the assumptions of Lemma~\ref{lem:compression}, the relative approximate error of the 
		reduced approximation satisfies
		\begin{equation}
			\frac{\|X^* - V\hat{X}V^\top\|_F}{\|V\hat{X}V^\top\|_F}
			\le \frac{\varepsilon(d_{\mathrm{comp}})}{\|\hat{X}\|_F},
			\label{eq:main_bound}
		\end{equation}
		with $\varepsilon(\cdot)$ as in \eqref{eq:epsilon}.
	\end{proposition}
	
	\begin{proof}
		From Lemma~\ref{lem:compression}, $\|\Delta\|_F \le \varepsilon(d)$. 
		Since $d \le d_{\mathrm{comp}}$ and $\varepsilon$ is increasing, 
		$\varepsilon(d) \le \varepsilon(d_{\mathrm{comp}})$. Dividing by 
		$\|V\hat{X}V^\top\|_F = \|\hat{X}\|_F$ gives \eqref{eq:main_bound}.
	\end{proof}
	
	\begin{remark}[Role in the basis construction]
		\label{rem:role_in_basis}
		The role of Proposition~\ref{thm:main} is to reveal the mechanism by which the choice of the subspace \(V\) controls the total error: since \(\varepsilon(\cdot)\) is increasing in \(d\), minimizing the aggregate training projection error, which in turn controls the indicator \(d_{\text{comp}}\) and hence \(\varepsilon(d_{\text{comp}})\), indirectly compresses the error bound for new parameters. This observation motivates the Grassmann manifold optimization in Section~\ref{sec:framework}, which is designed to minimize the projection error over the Grassmann manifold. The role of the residual-based certificate in Subsection~\ref{subsec:residual_estimation} is complementary: it furnishes an online bound once a reduced solution is available, whereas the analysis of this subsection explains \emph{why} the subspace chosen by Grassmann optimization yields a small value of that bound.
	\end{remark}

    To make the link between the Grassmann objective $\Phi(V)$ and the projection error indicator $d_{\mathrm{comp}}$ explicit, we establish the following simple bound. Define the snapshot norm scale
    \[
        M_X := \max_{1 \le j \le n_s} \|X^*(\mu_j)\|_F
    \]
    and the training-set fill distance
    \[
        h_{\mathrm{train}} := \sup_{\mu \in \mathcal{P}} \min_{1 \le j \le n_s}\|\mu - \mu_j\|_2.
    \]
    \begin{corollary}[Grassmann objective controls $d_{\mathrm{comp}}$]\label{cor:grassmann_to_dcom}
        For any $\mu \in \mathcal{P}$, the projection error indicator satisfies
        \begin{equation}
            d_{\mathrm{comp}} \le M_X \sqrt{\Phi(V)} + L_f \, h_{\mathrm{train}}.
            \label{eq:grassmann_link}
        \end{equation}
    \end{corollary}

    \begin{proof}
    By definition~\eqref{eq:proj_bound},
    $d_{\mathrm{comp}} = f(\mu_{\mathrm{nn}}) + L_f \|\mu - \mu_{\mathrm{nn}}\|_2$. For the first term, since $\Phi(V)$ is the sum of squared relative projection errors over all training snapshots, each individual summand is bounded by the total:
    \[
        \frac{\|X^*(\mu_{\mathrm{nn}}) - P X^*(\mu_{\mathrm{nn}}) P\|_F^2}
        {\|X^*(\mu_{\mathrm{nn}})\|_F^2} \le \Phi(V).
    \]
    Taking the square root and multiplying by
    $\|X^*(\mu_{\mathrm{nn}})\|_F \le M_X$ gives
    \[
        f(\mu_{\mathrm{nn}})
        = \|X^*(\mu_{\mathrm{nn}}) - P X^*(\mu_{\mathrm{nn}}) P\|_F
        \le M_X \sqrt{\Phi(V)}.
    \]
    For the second term, the definition of $h_{\mathrm{train}}$ directly gives
    $\|\mu - \mu_{\mathrm{nn}}\|_2 \le h_{\mathrm{train}}$. Combining both inequalities yields the result.
    \end{proof}
    This corollary clarifies the mechanism behind the basis construction: the Grassmann objective $\Phi(V)$ controls the training projection errors, which in turn controls $d_{\mathrm{comp}}$ through two terms. The first term $M_X \sqrt{\Phi(V)}$ is reduced by optimizing $V$ on the Grassmann manifold; the second term $L_f h_{\mathrm{train}}$ depends only on the training-set coverage and the Lipschitz constant, not on $V$. Since $\varepsilon(\cdot)$ is monotonically increasing, a smaller $d_{\mathrm{comp}}$ yields a correspondingly smaller $\varepsilon(d_{\mathrm{comp}})$ in Proposition~\ref{thm:main}, which quantifies how the subspace choice propagates into the overall approximation error. In practice, $L_f$ is replaced by the empirical estimate $\hat{L}_f$ from Remark~\ref{rem:Lf_computation}.

	\section{Numerical experiments}\label{sec:numerics}
	
	In this section, we evaluate the performance of the proposed global and localized reduced-order methods on two test cases. The first example, taken from the literature, serves as a benchmark for the global basis approach and illustrates the behavior of the error estimator on a standard problem. The second example is a synthetically constructed one-dimensional heat conduction problem designed to produce a clustered, multi-modal solution manifold, which highlights the advantages of the localized basis method. 

    The CPU time measurements reported in this work were performed using MATLAB on the High Performance Computing Platform of Xiangtan University.

    In all experiments, the Riemannian conjugate-gradient iterations on the Grassmann manifold are initialized with the one-sided optimal solution: the basis composed of the leading \(r\) eigenvectors of the matrix \(\sum_j w_jX^*(\mu_j)^2\)\cite{DingY05}.
	
	\subsection{Global Basis}
	
	This model is taken from \cite{Kressner}. Consider the parameterized heat equation
	\begin{equation}
		\begin{aligned}
			\frac{\partial \vartheta}{\partial t} - \nabla \cdot \bigl(\sigma(\xi) \nabla \vartheta\bigr) &= f \quad \text{in } \Omega \times (0,T),\\
			\vartheta &= 0 \quad \text{on } \partial\Omega \times (0,T),
		\end{aligned}
		\label{eq:heat_pde}
	\end{equation}
	on the spatial domain $\Omega = (0,4)^2 \subset \mathbb{R}^2$, with source term $f \equiv 1$. The thermal conductivity field $\sigma(\xi)$ is piecewise constant and takes the form
	\begin{equation}
		\sigma(\xi) =
		\begin{cases}
			1 + \mu_i, & \xi \in D_i, \quad i = 1,\dots,4,\\
			1, & \xi \in \Omega \setminus \bigl(\cup_{i=1}^4 D_i\bigr),
		\end{cases}
		\label{eq:conductivity}
	\end{equation}
	where $D_i \subset \Omega$, $i = 1,\dots,4$, are four disks of radius $0.5$ centered at $(1,1)$, $(3,1)$, $(1,3)$, and $(3,3)$, respectively. The parameter vector $\mu = [\mu_1,\mu_2,\mu_3,\mu_4]^\top$ varies in $\mathcal{P} = [0.1,2]^4$.
	
	Discretizing \eqref{eq:heat_pde} by the finite element method with piecewise linear basis functions yields a system of dimension $N = 1580$, which is then converted into the PGARE form \eqref{eq:PGARE}. The resulting mass matrix $E(\mu) \equiv E \in \mathbb{R}^{N \times N}$ is symmetric positive definite and parameter-independent. The stiffness matrix $A(\mu)$ admits the affine decomposition
	\begin{equation}
		A(\mu) = \mu_1 A_1 + \mu_2 A_2 + \mu_3 A_3 + \mu_4 A_4 + A_5,
		\label{eq:stiffness_affine}
	\end{equation}
	where $A_i \in \mathbb{R}^{N \times N}$, $i = 1,\dots,4$, are symmetric negative semidefinite and $A_5$ is symmetric negative definite. The coefficient matrices $B(\mu)$ and $C(\mu)$ are constructed to satisfy $G(\mu) = B(\mu)B(\mu)^\top \succeq 0$ and $Q(\mu) = C(\mu)^\top C(\mu) \succeq 0$.
	
\begin{figure}
		\centering
		\includegraphics[width=0.8\linewidth]{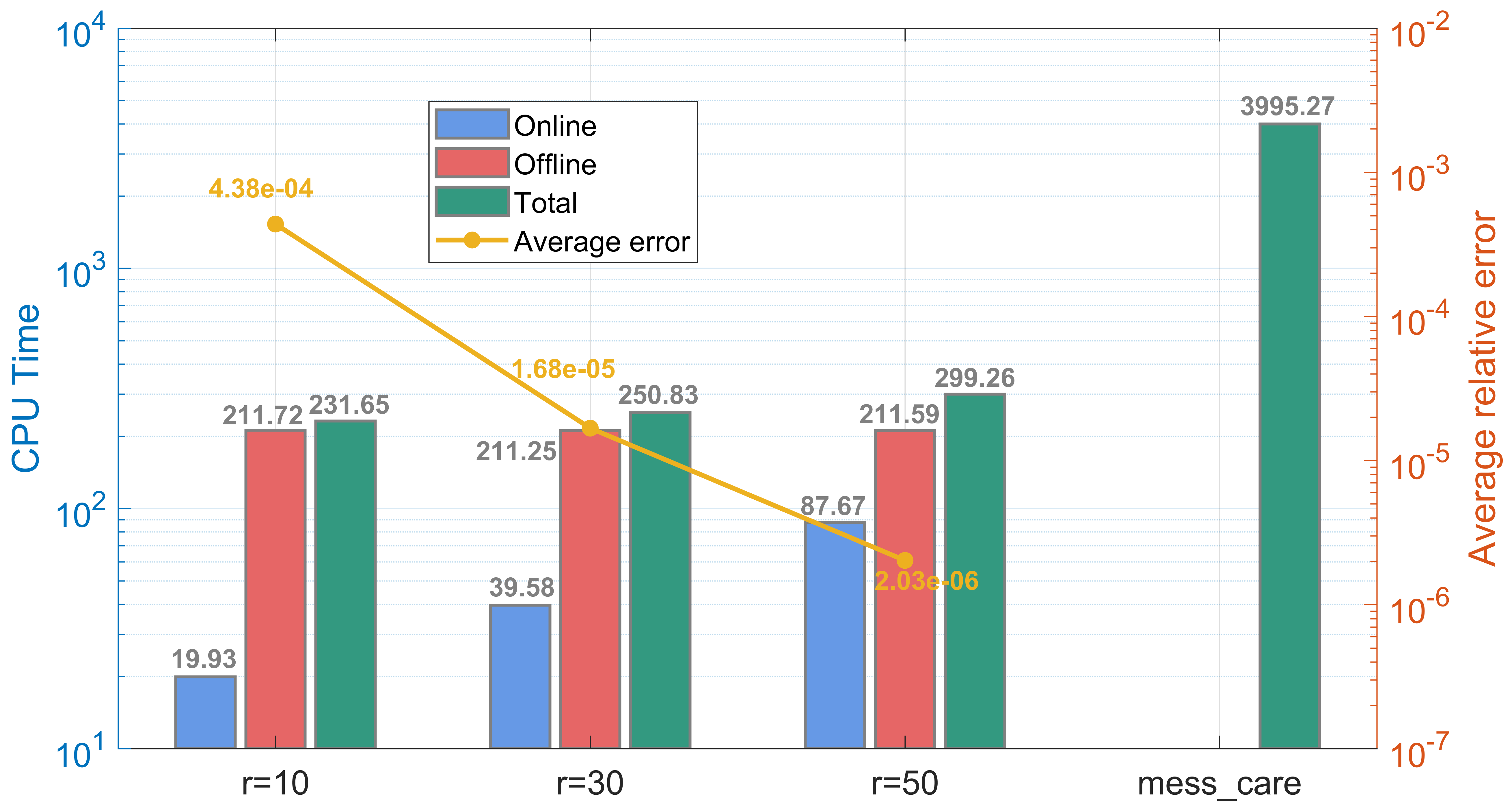}
		\caption{Average relative errors of the proposed method with \(r=10,30,50\), and CPU time comparison with the full-order solver \text{mess\_care} of the M-M.E.S.S. toolbox, evaluated over \(10^4\) test parameters.}
		\label{fig:global_perf}
	\end{figure}

    We first examine the performance of the global reduced basis approach. Fig.~\ref{fig:global_perf} reports the average relative errors of the proposed method with reduced orders \(r=10,30,50\), alongside a CPU time comparison with the full-order solver \text{mess\_care} of the M-M.E.S.S. toolbox \cite{saak2023mmess}. The reduced bases are constructed using a training set of 500 parameters, and the performance is evaluated over \(10^4\) test parameters. As expected, increasing the reduced dimension improves the accuracy of the reduced-order solution. Notably, the proposed method achieves a good speedup compared to the full-order solver, demonstrating the computational efficiency of the reduced-order framework for global bases.

\begin{figure}
		\centering
		\begin{subfigure}{0.49\textwidth} 
			\includegraphics[width=\linewidth]{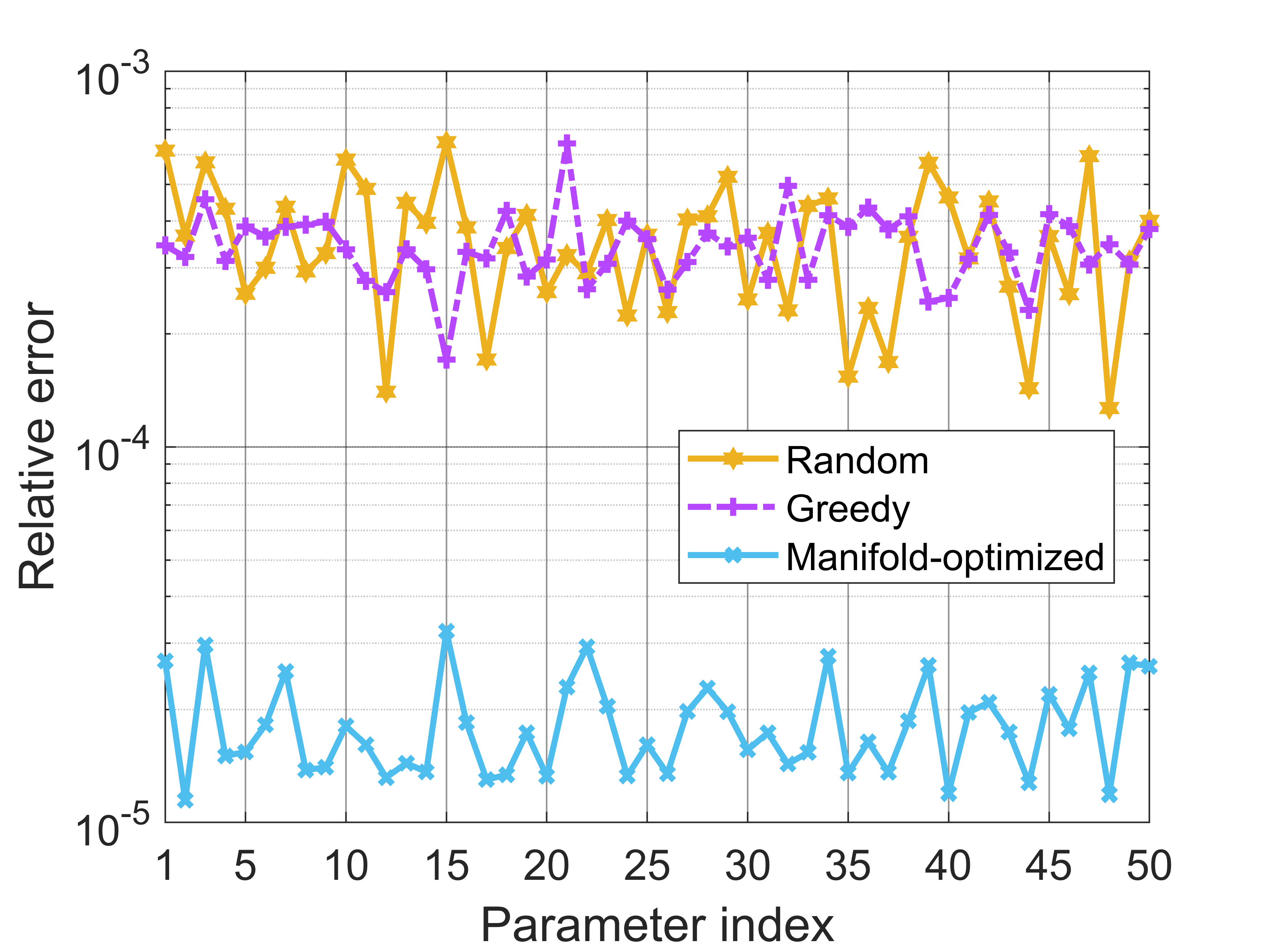}
			\caption{Pointwise error}
		\end{subfigure}
		\hspace{-0.1cm}
		\begin{subfigure}{0.49\textwidth} 
			\includegraphics[width=\linewidth]{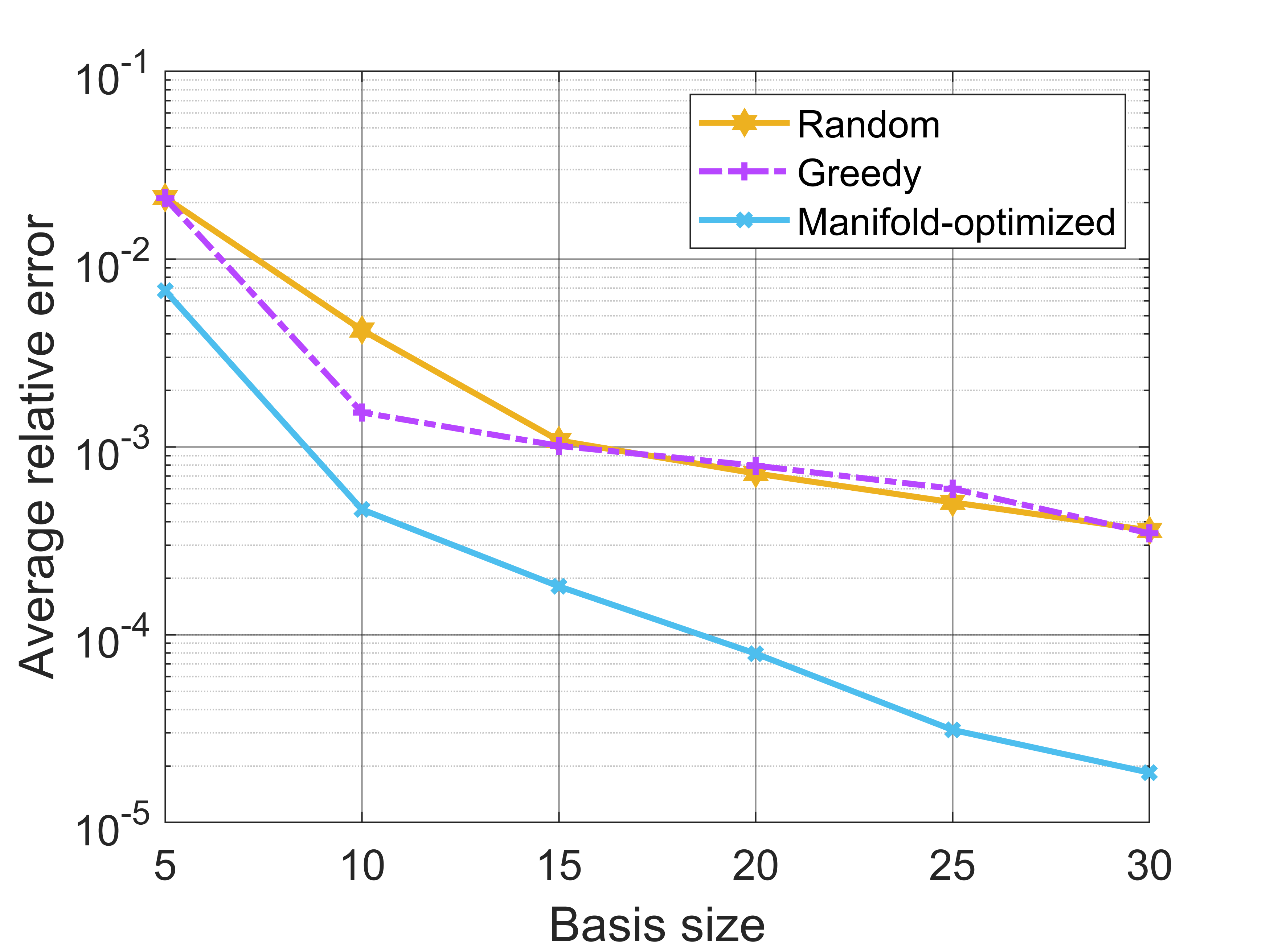}
			\caption{Average relative error}
		\end{subfigure}
		\caption{Comparison of the random, greedy, and manifold-optimized bases: (a) Pointwise relative errors at 50 test parameters for reduced orders $r=30$; (b) Average relative errors as functions of the reduced order $r$, evaluated over the test parameters.}
		\label{Greedy_manifold_Random}
	\end{figure}

    Fig.~\ref{Greedy_manifold_Random} compares the three basis-construction strategies on the same training snapshots. For simplicity, the greedy strategy used here iteratively selects the snapshot with the largest exact error (computed directly rather than via an error estimator), whereas the random approach samples an equal number of snapshots randomly. At the fixed reduced order \(r=30\), the manifold-optimized basis attains pointwise relative errors on the order of \(10^{-5}\) across all 50 test parameters, about one order of magnitude below the random and greedy bases, whose errors fluctuate near \(10^{-4}\)--\(10^{-3}\). As \(r\) increases from \(5\) to \(30\), all three methods reduce their average errors, but the manifold-optimized curve consistently lies well below the other two and decreases more rapidly, reaching \(2\times10^{-5}\) at \(r=30\). These results indicate that optimizing the aggregate two-sided projection error directly on the Grassmann manifold yields a markedly more accurate reduced space than either random or greedy sampling.
	
\begin{figure}
		\centering
		\begin{subfigure}{0.49\textwidth} 
			\includegraphics[width=\linewidth]{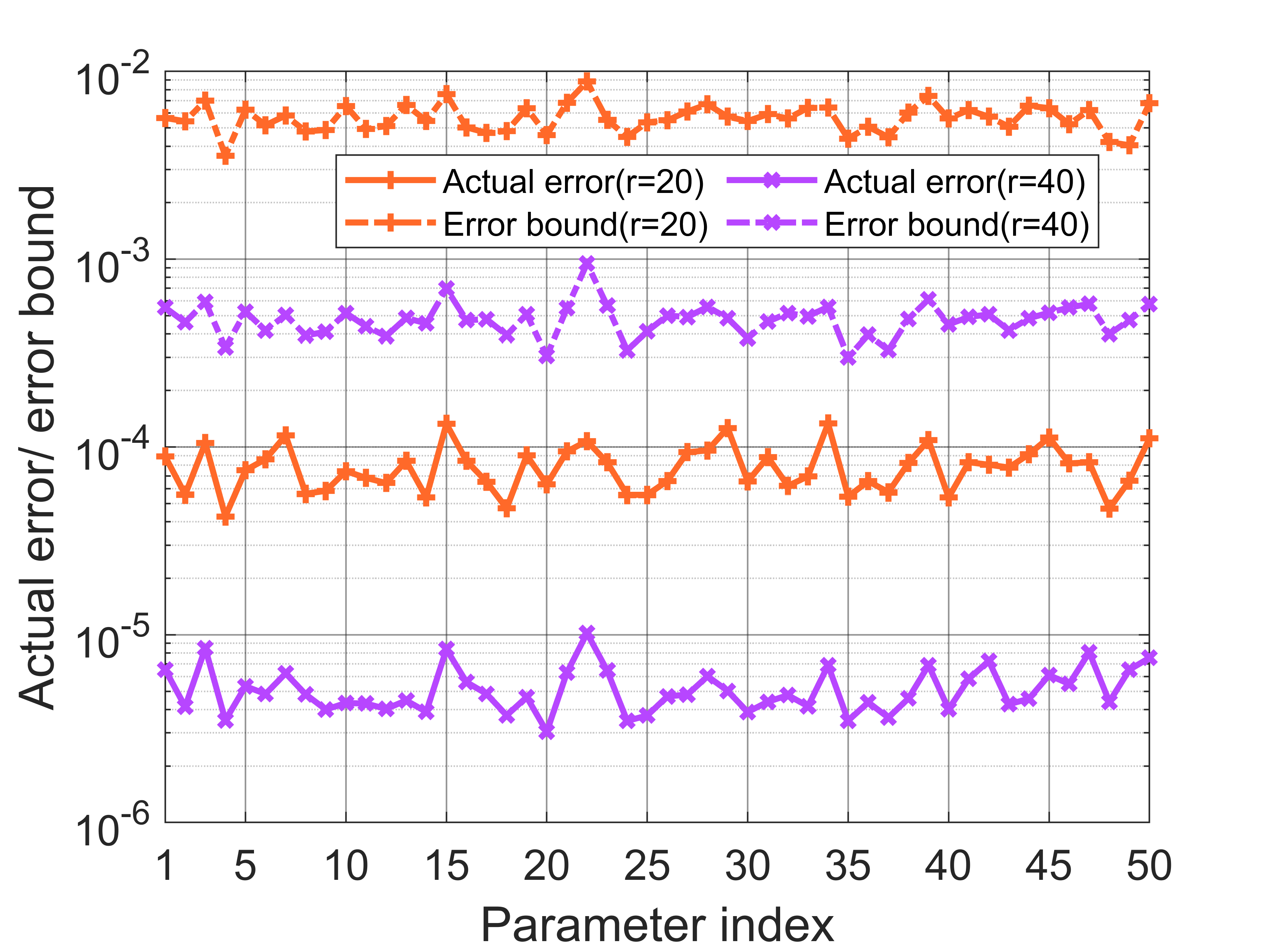}
			\caption{Pointwise error and error bound}
		\end{subfigure}
		\hspace{-0.5cm}
		\begin{subfigure}{0.49\textwidth} 
			\includegraphics[width=\linewidth]{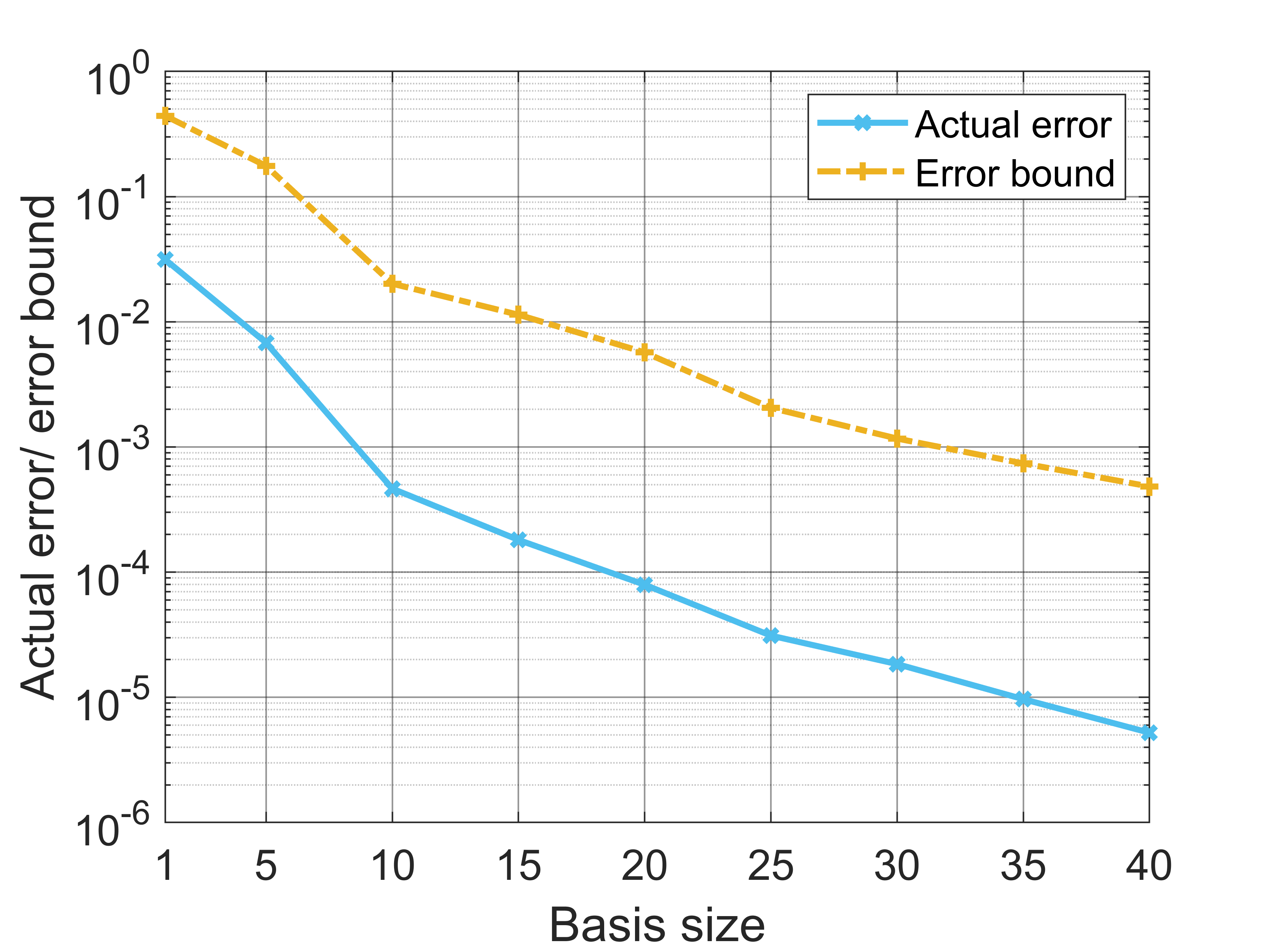}
			\caption{Actual error and error bound}
		\end{subfigure}
		\caption{(a) Relative errors and error bounds at test parameters for reduced orders $r=20,40$; (b) Average relative errors and their bounds as functions of the reduced order $r$, evaluated over the test parameters.}
		\label{20_40}
	\end{figure}

    Fig.~\ref{20_40} presents two complementary views. Panel (a) shows the pointwise relative errors for reduced orders \(r=20\) and \(r=40\) across the 50 test parameters, together with their corresponding error estimator. The error bounds consistently lie above the actual errors for all certified points, confirming the validity of the theoretical bound. Panel (b) displays the average relative errors and their bounds as functions of the reduced order \(r\), averaged over the same 50 test parameters. As \(r\) increases, both the actual errors and the bounds decrease, with the bounds remaining consistently conservative.
	
\begin{figure}
		\centering
		\begin{subfigure}{0.49\textwidth} 
			\includegraphics[width=\linewidth]{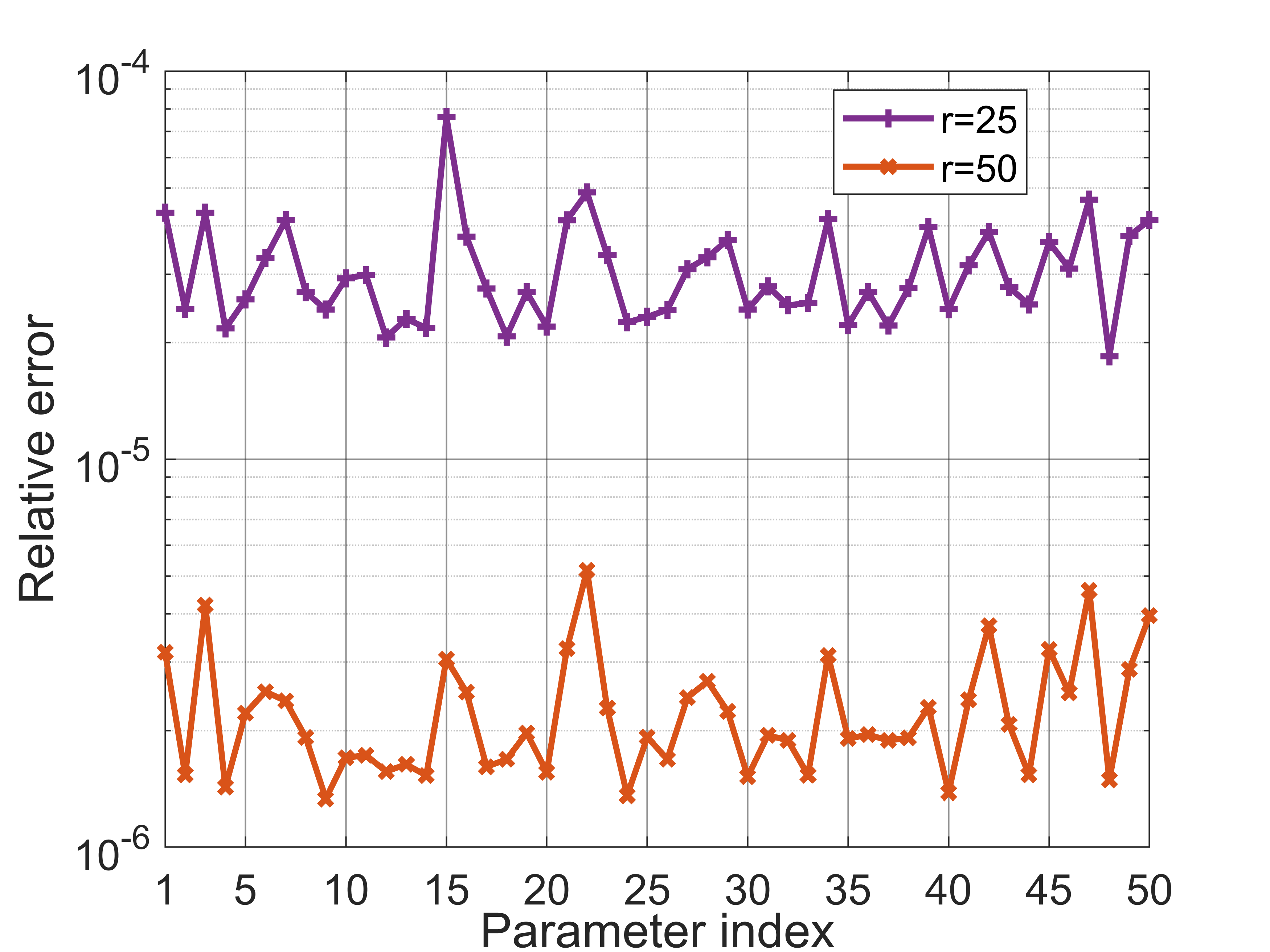}
			\caption{Actual pointwise relative errors}
		\end{subfigure}
		\hspace{-0.5cm}
		\begin{subfigure}{0.49\textwidth} 
			\includegraphics[width=\linewidth]{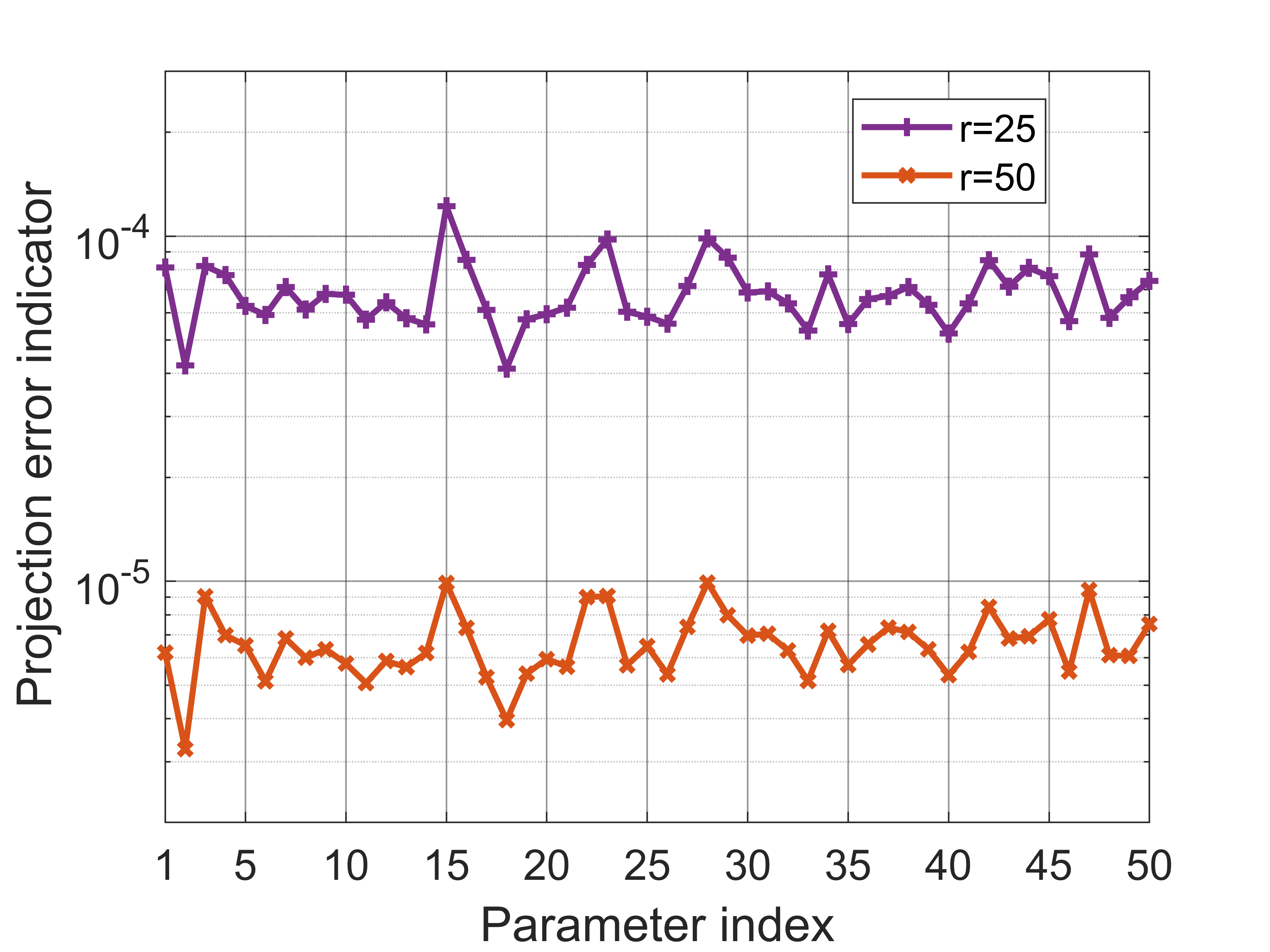}
			\caption{Projection error indicator}
		\end{subfigure}
		\caption{(a) Actual relative errors at test parameters for reduced orders $r=25,50$; (b) projection error indicator at test parameters for reduced orders $r=25,50$.}
		\label{25_50}
	\end{figure}

	Fig.~\ref{25_50} displays the actual relative errors and the projection error indicator over the 50 test parameters for reduced orders $r=25, 50$. Panel (a) reports the actual relative errors, while panel (b) shows the corresponding projection error indicators. As the reduced dimension increases from 
	25 to 50, both the actual relative error and the projection error indicator decrease noticeably across all test parameters. The consistent reduction in both quantities indicates that the manifold-optimized basis construction effectively captures the dominant features of the Riccati solution manifold, thereby validating the effectiveness of the proposed Grassmann optimization approach.

	
	
	\subsection{Local Basis}
	
	We now consider a one-dimensional heat conduction problem that is specifically designed to produce a solution manifold with a clustered, multi-modal geometry. Unlike the previous example, which is taken from the literature, this test case is constructed synthetically to clearly demonstrate the advantage of localized reduced bases when the solution manifold exhibits multiple separated clusters.

    Consider the one-dimensional heat conduction equation on the spatial domain $\Omega = (0, L)$ with $L = 1$:
    \begin{equation}\label{eq:heat-pde}
	    \frac{\partial \vartheta}{\partial t} - \frac{\partial}{\partial x}\!\left(\sigma(x;\mu)\, \frac{\partial \vartheta}{\partial x}\right) = 0 \quad \text{in } (0,L)\times(0,T),
    \end{equation}
    subject to homogeneous Neumann (zero-flux) boundary conditions
    \begin{equation}\label{eq:heat-bc}
	    \sigma(0;\mu)\frac{\partial \vartheta}{\partial x}(0, t) = \sigma(L;\mu)\frac{\partial \vartheta}{\partial x}(L, t) = 0.
    \end{equation}
    The thermal conductivity field $\sigma(x;\mu)$ is piecewise constant and admits the affine decomposition
    \begin{equation}\label{eq:heat-sigma}
	    \sigma(x;\mu) = \chi_{R_1}(x) + \mu_1\, \chi_{R_2}(x) + \mu_2\, \chi_{R_3}(x) + \mu_3\, \chi_{R_4}(x) + \mu_4\, \chi_{R_5}(x),
    \end{equation}
    where $\{R_i\}_{i=1}^5$ are five equal subintervals partitioning $(0,1)$ and $\chi_{R_i}$ denotes the characteristic function of $R_i$. The conductivity of region~$R_1$ is fixed at $1$, while the conductivities of regions $R_2$--$R_5$ are controlled by the parameters $\mu_1, \ldots, \mu_4 \in [0,1]$. The inclusion of the boundary value $\mu_i = 0$ is deliberately designed to model perfectly insulated subregions. Even in these degenerate cases, the non-vanishing conductivity in $R_1$ together with the damping term $-\gamma_d I_N$ in \eqref{eq:heat-A} ensures that the system matrix $A(\mu)$ remains strictly dissipative, thereby guaranteeing the well-posedness of the PGARE.

    The control input $u(t)$ enters through a spatially distributed heat source, and the output is a pointwise temperature measurement. The control and observation matrices are constructed as convex combinations of five fixed profiles:
    \begin{equation}\label{eq:heat-BC}
	    B(\mu) = \sum_{i=1}^5 \frac{\alpha_i(\mu)}{10}B_i, \qquad C(\mu) = \sum_{i=1}^5 \alpha_i(\mu)\, C_i,
    \end{equation}
    where $B_i \in \R^{N\times 1}$ are discretized Gaussian source profiles with standard deviation $\sigma_B=0.5$, centered at five equally spaced locations in $(0,1)$, $C_i \in \R^{1\times N}$ are point-evaluation functionals at the same locations, and the weights $\alpha_i(\mu)$ are normalized radial basis functions:
    \begin{equation}\label{eq:heat-rbf}
	    \alpha_i(\mu) = \frac{\exp\!\bigl(-\|\mu - c_i\|^2 / (2\gamma_{\mathrm{rbf}}^2)\bigr)}{\sum_{j=1}^5 \exp\!\bigl(-\|\mu - c_j\|^2 / (2\gamma_{\mathrm{rbf}}^2)\bigr)},\ \gamma_{\mathrm{rbf}} = 0.1.
    \end{equation}
    The five parameter-space centers
    \begin{equation}\label{eq:heat-centers}
	    c_1 = \begin{bmatrix}1\\0\\0\\0\\0.2\end{bmatrix},\;
	    c_2 = \begin{bmatrix}0\\1\\0\\0\\0.4\end{bmatrix},\;
	    c_3 = \begin{bmatrix}0\\0\\1\\0\\0.6\end{bmatrix},\;
	    c_4 = \begin{bmatrix}0\\0\\0\\1\\0.8\end{bmatrix},\;
	    c_5 = \begin{bmatrix}0\\0\\0\\0\\1.0\end{bmatrix}
\end{equation}
    represent five distinct operating regimes. The first four centers raise the conductivity of the subregion $R_{i+1}$ to $\mu_i = 1$, while the fifth center leaves all four regional conductivities at their lower bound, so that the conductivity field reduces to the baseline $\sigma = \chi_{R_1}$, and instead shifts the actuation/sensing profile through $\mu_5$. The well-separated centers naturally induce a clustered, multi-modal solution manifold, which is the target scenario for the localized basis method.

    The parameter domain is \(\mathcal{P} = [0,1]^5\). For this synthetic benchmark, we construct the semi-discrete system directly in the standard state-space form with \(E = I_N\). The spatial domain is discretized by uniform piecewise linear finite elements with mesh size \(h = 1/1200\), which yields \(N = 1201\) degrees of freedom; the subregion interfaces $x = 0.2, 0.4, 0.6, 0.8$ fall on element boundaries, so that the region-wise stiffness matrices are assembled exactly. In the finite element discretization, the zero-flux condition \eqref{eq:heat-bc} is naturally satisfied by omitting the boundary terms in the weak formulation, which is consistent with the assembled stiffness matrices $K_i$ used in our implementation. The stiffness matrices \(K_i\) are assembled from the finite element discretization of the diffusion operator on the five subregions (with unit conductivity) and are used directly to assemble the state matrix. The system matrix admits the affine decomposition
    \begin{equation}\label{eq:heat-A}
	    A(\mu) = -\bigl(K_0 + \mu_1 K_1 + \mu_2 K_2 + \mu_3 K_3 + \mu_4 K_4\bigr) - \gamma_d I_N,
    \end{equation}
    where $\gamma_d = 0.2$. The first four parameter components $\mu_1, \dots, \mu_4$ control the regional conductivities in $R_2, \dots, R_5$, whereas the fifth component $\mu_5$ does not enter $A(\mu)$ directly but affects the actuator/observation profile $B(\mu), C(\mu)$ through the RBF weights $\alpha_i(\mu)$. The term $-\gamma_d I_N$ ensures that all eigenvalues of $A(\mu)$ have negative real part for every $\mu \in \mathcal{P}$, guaranteeing the existence of a unique stabilizing solution $X^*(\mu)$ to the PGARE.

    Training and test parameters are sampled around the five parameter-space regime centers~\eqref{eq:heat-centers} with Gaussian noise:
    \begin{equation}\label{eq:heat-sampling}
	    \mu = c_i + \varepsilon, \qquad \varepsilon \sim \mathcal{N}(0, \sigma_{\mathrm{sample}}^2\, I_5), \qquad \sigma_{\mathrm{sample}} = 0.05,
    \end{equation}
    clamped to $[0,1]^5$, and labeled by the generating center $c_i$. The training set consists of $100$ parameters per center, i.e., $n_s = 500$ parameters in total, and every method compared below is computed on the same full-order snapshots $X^*(\mu_j)$, $j = 1, \dots, n_s$.

    For the localized reduction we use $n_c = 5$ clusters, a local basis dimension $r_\ell = 30$, and a hybrid clustering weight $\alpha = 0.5$, this fixed value balancing the parameter-space distance and the solution-space correlation in our problem setting; the $k$-medoids algorithm is run with 200 random initializations. For online model selection in this synthetic benchmark, we employ a linear discriminant analysis (LDA) classifier \cite{Hastie2009}, trained on the same $n_s = 500$ labeled parameters, whose classification accuracy is assessed on the independent test parameters. Since the training parameters are sampled as Gaussian clouds around the five centers with a common covariance structure, the assumptions underlying LDA are approximately satisfied, making it a natural and efficient choice for this specific example. We note that the localized framework itself is independent of the choice of classifier; other classifiers could be used in situations where the cluster boundaries are more complex.

\begin{figure}
		\centering
		\begin{subfigure}{0.55\textwidth} 
			\includegraphics[width=\linewidth]{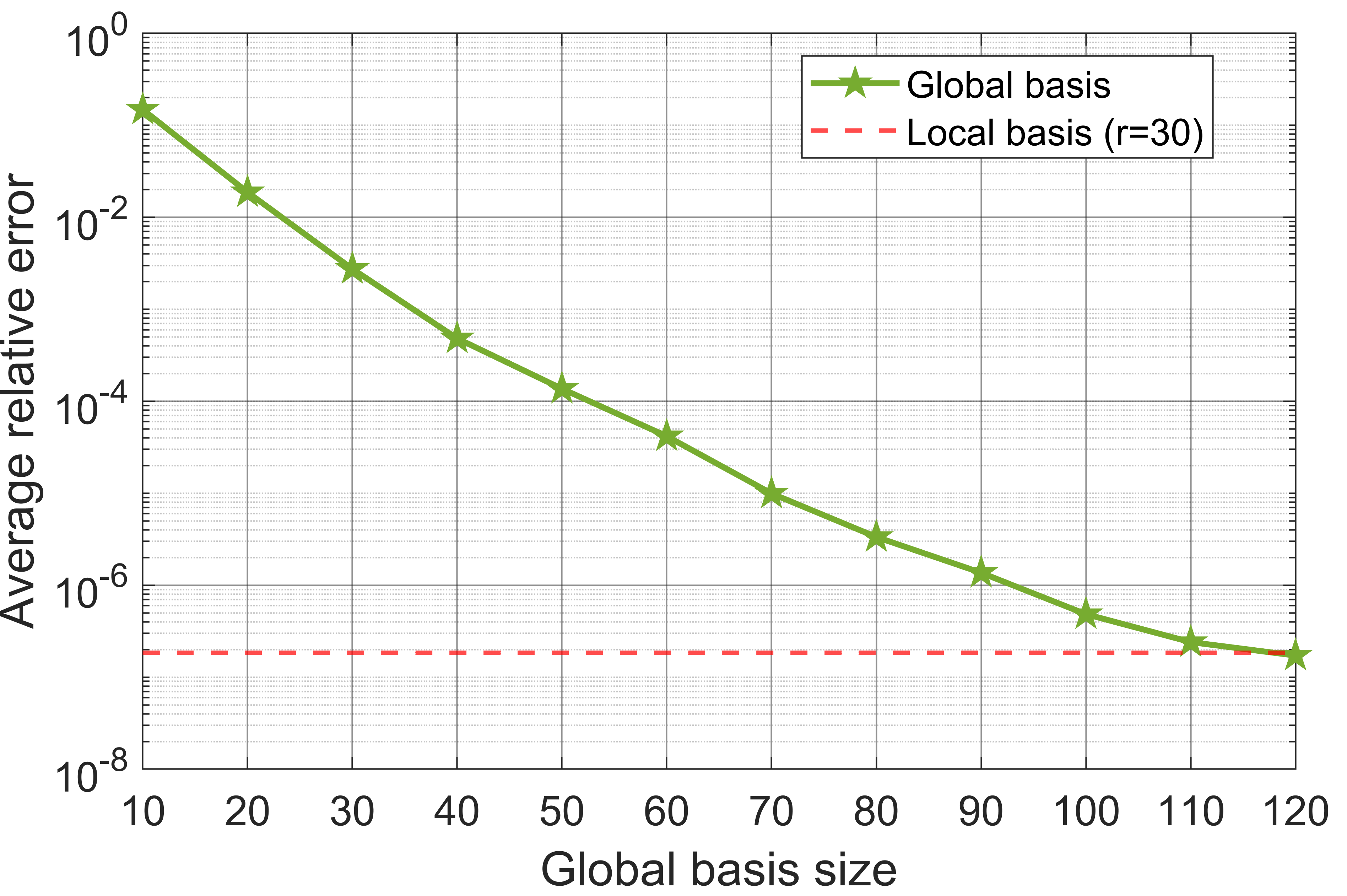}
			\caption{Average relative errors}
		\end{subfigure}
		\hspace{-0.2cm}
		\begin{subfigure}{0.35\textwidth} 
			\includegraphics[width=\linewidth]{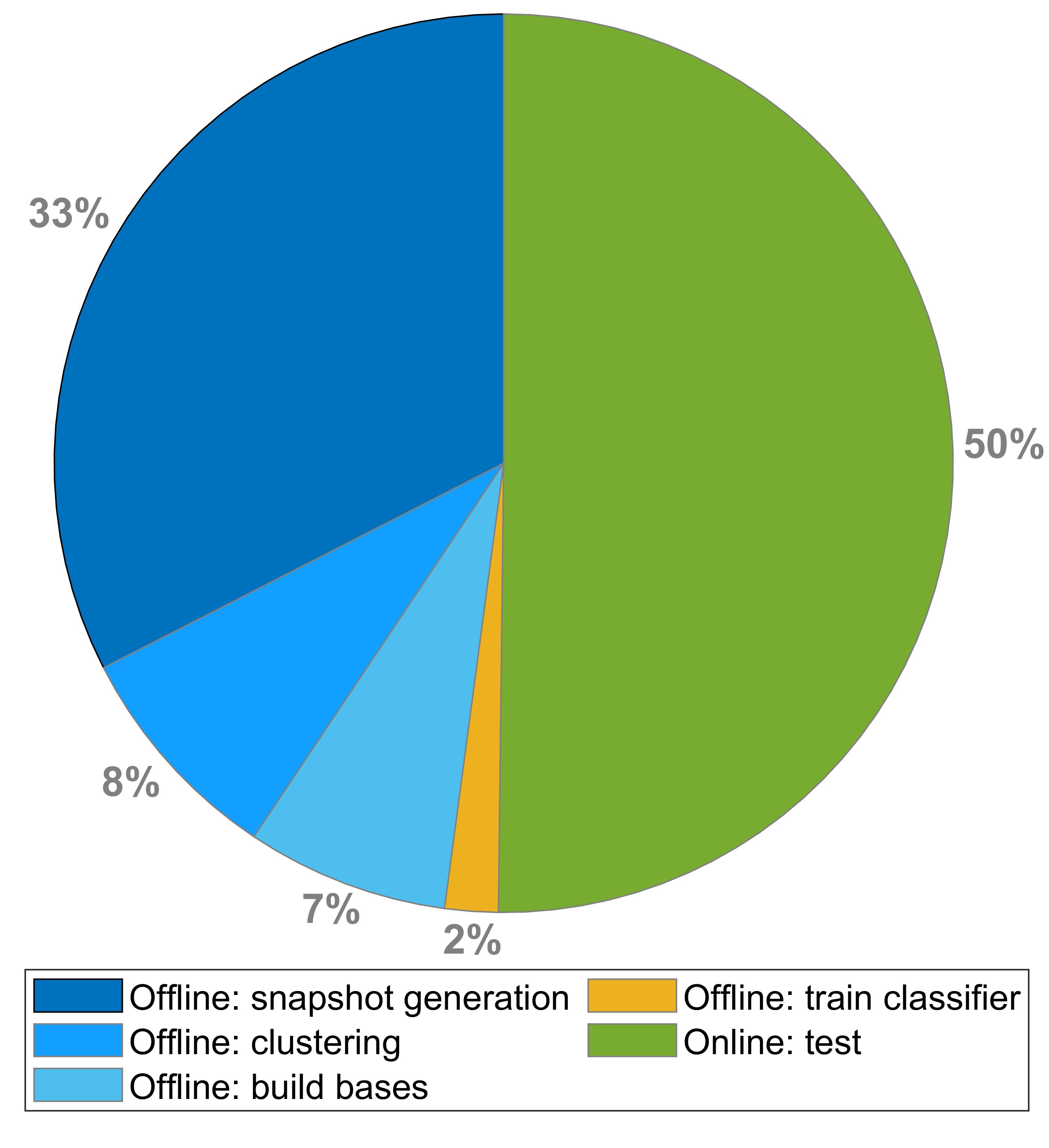}
			\caption{CPU time}
		\end{subfigure}
		\caption{(a) Average relative errors for local basis $(r=30)$, and global bases $(r=10 \ \text{to}\ 120)$, evaluated  over $10^4$ test parameters; (b) CPU time of offline and online stages, evaluated over $10^4$ test parameters.}
		\label{fig:local_vs_global}
	\end{figure}

	We now compare the localized basis approach with the global basis on this clustered problem. Fig.~\ref{fig:local_vs_global} presents the average relative errors and CPU time. Panel (a) shows that the local basis with \(r_{\text{local}}=30\) achieves an average relative error comparable to that of a global basis with \(r_{\text{global}}=120\), while requiring a substantially smaller online dimension per query. We note that the total number of stored basis vectors is $5 \times 30 = 150$, which exceeds the global dimension \(r_{\text{global}}=120\); the advantage of the localized approach lies in the reduced per-query online cost, not in the total storage. This demonstrates the efficiency of the localized approach: it attains high accuracy with a much smaller reduced dimension by exploiting the clustered structure of the solution manifold. Panel (b) reports the CPU time breakdown. The offline stage comprises the generation of snapshot solutions, clustering of training parameters, construction of local bases, and training of the classifier. The online stage consists of localizing the 10000 test parameters to their respective clusters and solving the reduced-order models; although this stage accounts for about half of the total CPU time, the cost per test parameter remains substantially lower than that of the full-order solver. 
	
\begin{figure}
		\centering
		\begin{subfigure}{0.49\textwidth} 
			\includegraphics[width=\linewidth]{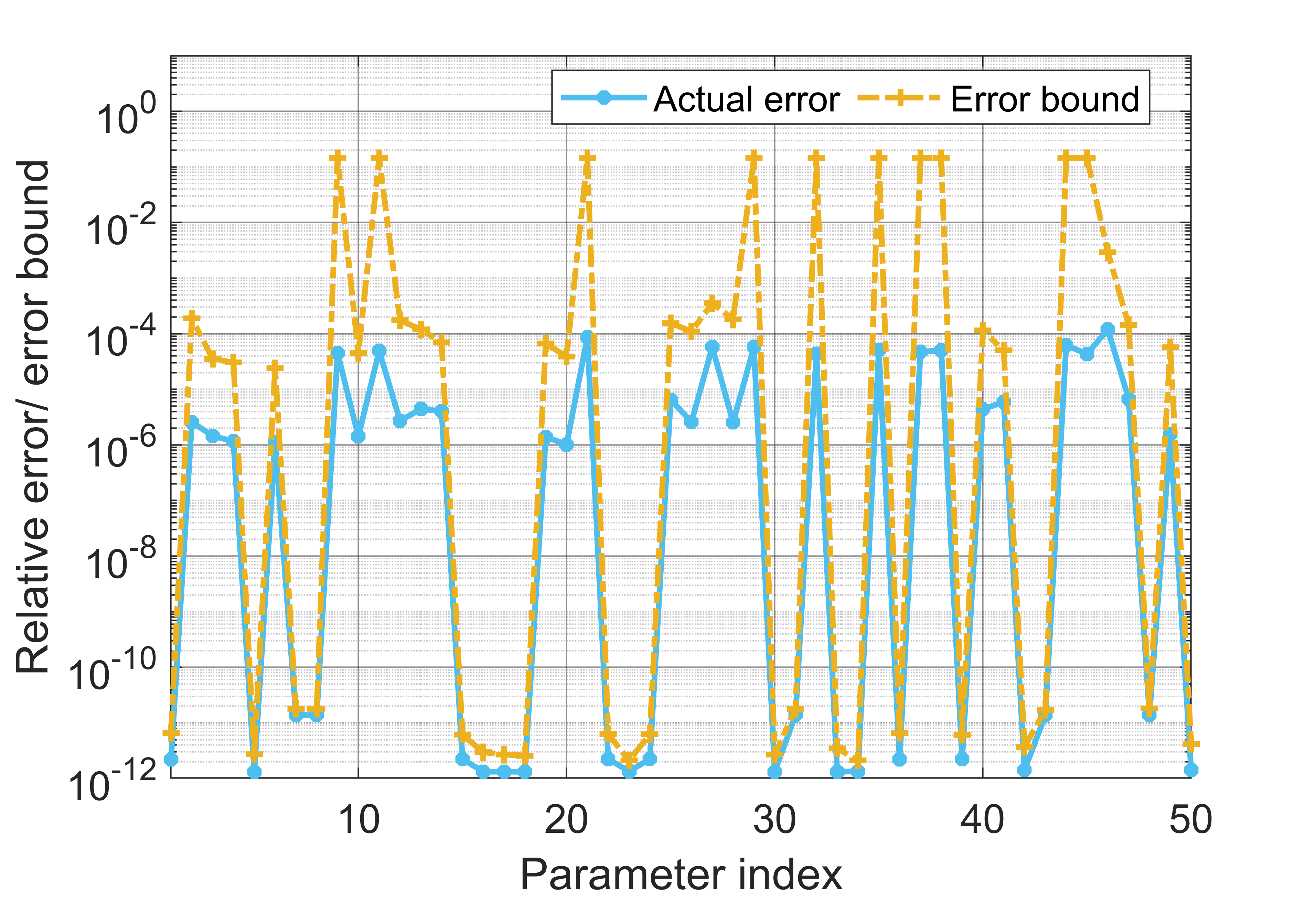}
			\caption{$r=15$}
		\end{subfigure}
		\hspace{-0.5cm}
		\begin{subfigure}{0.49\textwidth} 
			\includegraphics[width=\linewidth]{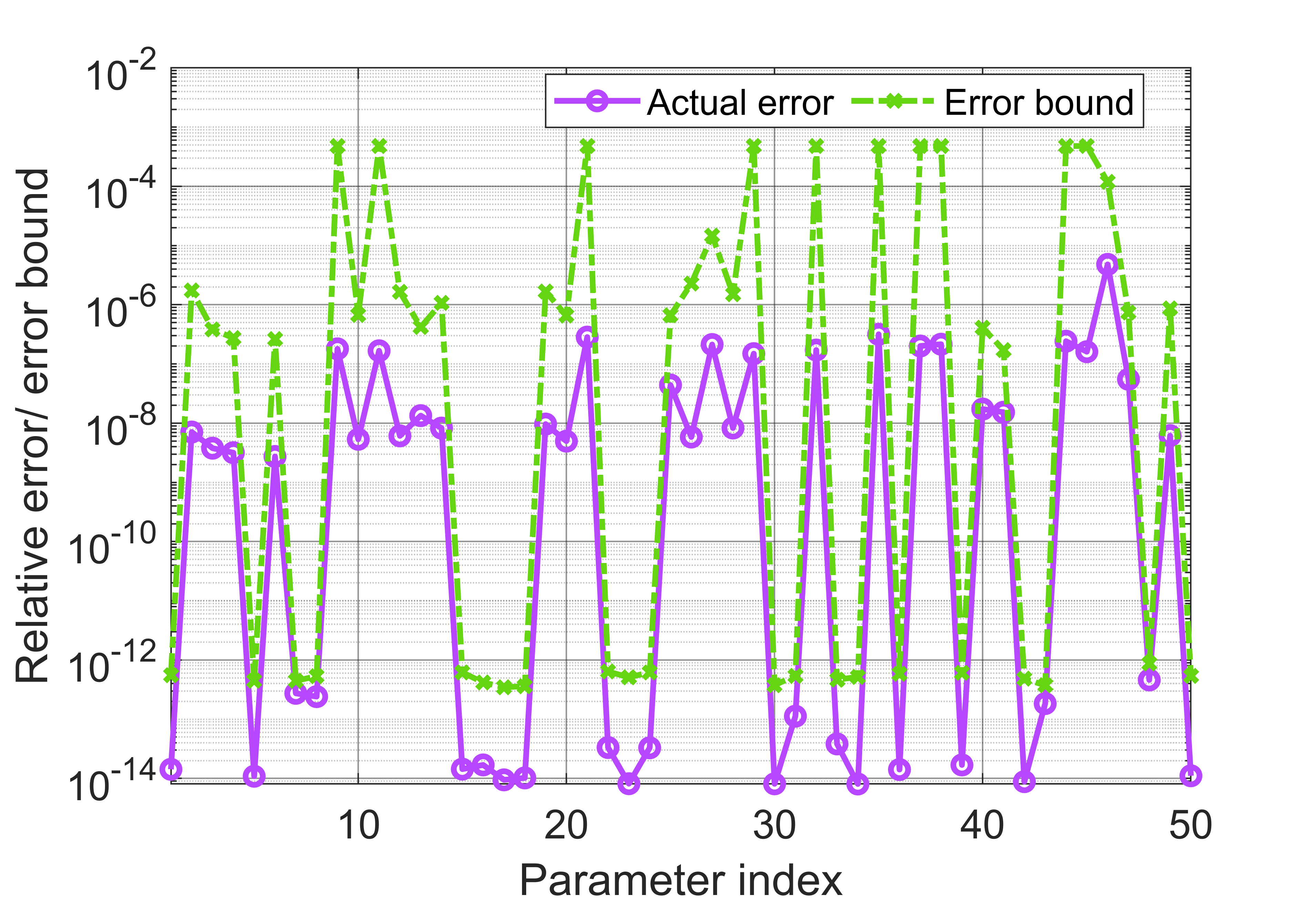}
			\caption{$r=30$}
		\end{subfigure}
		\caption{Actual error and error bound, evaluated over 50 test parameters.}
		\label{fig:local_bounds}
	\end{figure}

	The quality of the error bounds for the local basis is examined in Fig.~\ref{fig:local_bounds}. Panel (a) and (b) display the actual errors and the error bounds for the local basis with dimensions \(r_{\text{local}}=15, 30\) respectively. As expected, the bounds decrease as the reduced dimension increases, reflecting the improved approximation capability of larger subspaces. The error bounds consistently lie above the actual errors for all certified test points, confirming the reliability of the theoretical bound in the local setting.
	
\begin{figure}
		\centering
		\includegraphics[width=0.7\linewidth]{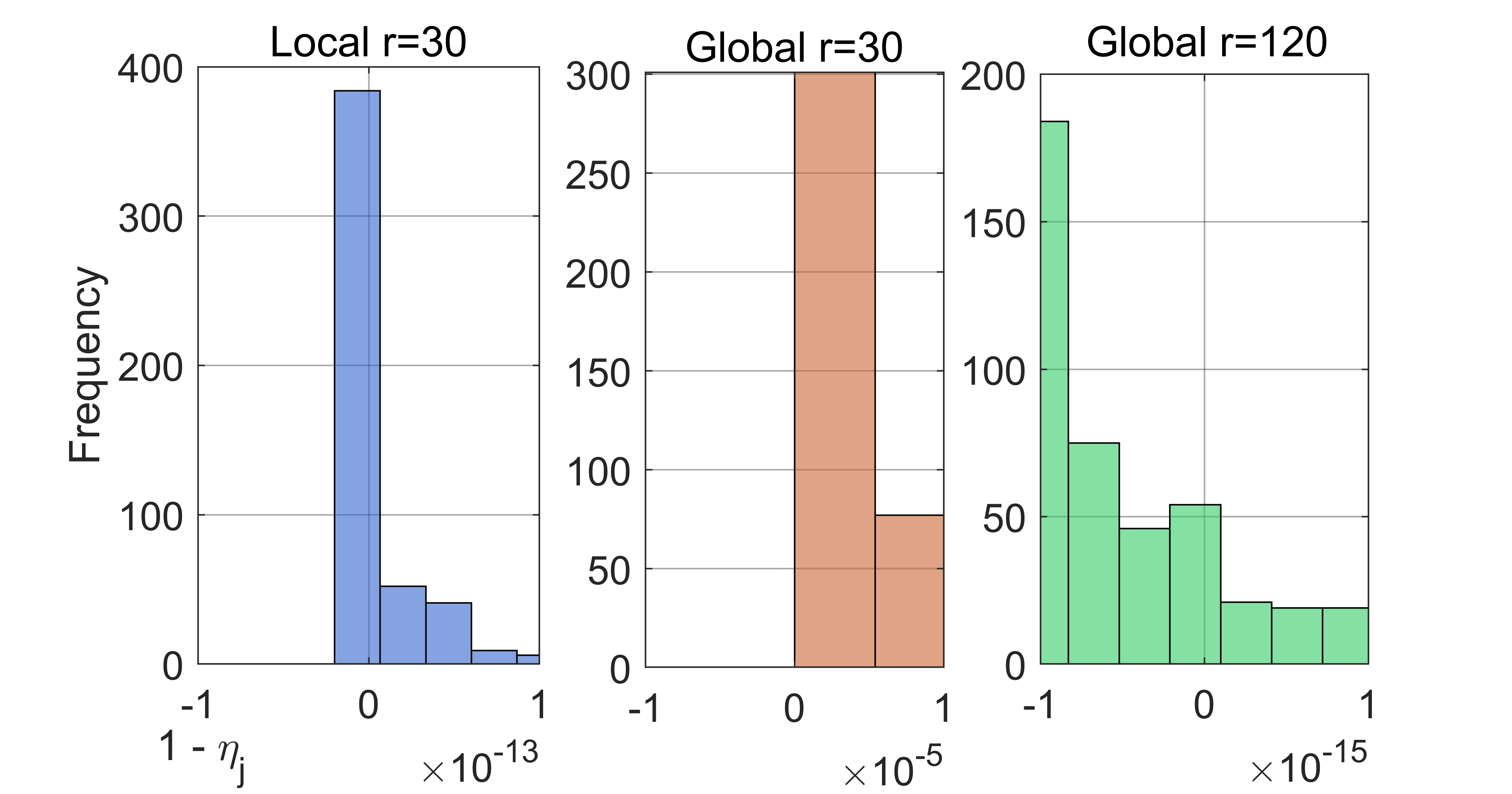}
		\caption{Uncaptured energy ratio $1-\eta_j$ at training parameters, for the local basis with $r_{\mathrm{local}}=30$, the global basis with $r_{\mathrm{global}}=30$, and the global basis with $r_{\mathrm{global}}=120$.}
		\label{fig:eta}
	\end{figure}

	To further assess the quality of the local subspaces, Fig.~\ref{fig:eta} plots the uncaptured energy ratio $1-\eta_j$ at each training point, where the energy capture ratio is defined as
    \[ 
        \eta_j := \frac{\|V^\top X^*(\mu_j) V\|_F^2}{\|X^*(\mu_j)\|_F^2},\ j = 1, \dots, n_s.
    \]
    By the Pythagorean identity \eqref{eq:pythagorean_identity}, the uncaptured ratio is equivalently
    \[
        1-\eta_j = \frac{\|X^*(\mu_j) - P X^*(\mu_j) P\|_F^2}{\|X^*(\mu_j)\|_F^2}.
    \]
    Results are shown for three different bases: the local basis with $r_{\text{local}}=30$, the global basis with $r_{\text{global}}=30$, and the global basis with $r_{\text{global}}=120$. The global basis with $r_{\text{global}}=30$ yields markedly larger $1-\eta_j$ than both the local basis with $r_{\text{local}}=30$ and the global basis with $r_{\text{global}}=120$, indicating that it leaves a substantially larger fraction of energy uncaptured. These results indicate that the local basis achieves energy capture comparable to the much larger global basis ($r_{\text{global}}=120$) at a substantially smaller active online dimension ($r_{\text{local}}=30$). This confirms that the clustering strategy effectively identifies and constructs subspaces that are well adapted to the local structure of the solution manifold, thereby reducing the online computational cost while maintaining accuracy.
	
	Finally, Table~\ref{tab:global_local_comparison} provides a quantitative summary of the performance comparison between the global and local bases over \(10^4\) test parameters. The local basis with \(r_{\text{local}}=30\) achieves an average relative error of \(1.85\times 10^{-7}\), which is comparable to the global basis with \(r=120\), while using a substantially smaller active online reduced dimension, at the cost of storing five cluster bases of dimension \(r_{\text{local}}=30\). The online test time for the local basis (37.73 seconds) is substantially lower than that of the global basis with \(r_{\text{global}}=120\) (502.35 seconds). These results indicate that, for clustered solution manifolds, the localized basis offers a more favorable accuracy–cost trade-off than a single global basis.
	
	\begin{table}
		\centering\small
		\caption{Performance comparison between global and local bases, evaluated over $10^4$ test parameters.}
		\label{tab:global_local_comparison}
		\begin{tabular}{c c c c c c}
			\toprule
			Basis &  Active dim & Total dim & Avg. rel. error & Train time (s) &Test time (s) \\
			\midrule
			\multirow{3}{*}{Global} & 30  & 30 & $2.75 \times 10^{-3}$ & 29.93 & 31.24 \\
			& 80 & 80  & $3.34 \times 10^{-6}$ & 28.86 & 191.16 \\
			& 120 & 120 & $1.72 \times 10^{-7}$ & 29.21 & 502.35 \\
			
			\midrule
			Local & 30 & 150  & $1.85 \times 10^{-7}$ & 37.43 & 37.73 \\
			\bottomrule
		\end{tabular}
	\end{table}
	
	\section{Conclusions}\label{sec:conclusion}
	In this paper, within the classical projection-based reduced-order modeling framework, we propose a method for constructing the projection basis via optimization on the Grassmann manifold for parameter-dependent generalized algebraic Riccati equations. The basis is obtained by solving a nonconvex optimization problem on the manifold, which preserves the intrinsic algebraic structure of the original system. Under the dissipativity assumptions of Lemma 2.1, the reduced equation admits a unique symmetric positive semidefinite stabilizing solution. The resulting affine reduced-order model achieves online computational complexity independent of the full-order dimension, making it well suited for many-query scenarios. To certify the reduced approximation, we recall the classical residual-based a posteriori error estimate from the literature. Furthermore, we decompose the approximation error to explain the rationale behind our basis construction. To handle clustered solution manifolds, we introduce a localized extension that clusters the solution snapshots to construct dedicated local bases. An online classifier then dynamically selects the appropriate local model for new parameters. Numerical experiments demonstrate that the framework achieves high accuracy and computational efficiency compared to the full-order solver. The localized models further reduce the active online dimension.

    Furthermore, the current localized framework does not account for potential classification errors introduced during the online stage. Since a misclassified parameter may lead to a suboptimal local basis and degrade accuracy, future work will focus on developing robust model-selection strategies and incorporating classification uncertainty into the online error estimation.
	
	\appendix
	
	\section{Proof of Theorem~\ref{thm:lipschitz}}
	\label{app:proof_lipschitz}
	
	\begin{proof}[Proof of Theorem~\ref{thm:lipschitz}]
			The proof consists of six steps.
			
			\paragraph{Step 1: Transformation to standard CARE} 
			For each $\mu$, define 
			\[
			Y(\mu):=E(\mu)^\top X(\mu)E(\mu).
			\]
			Substitution into the PGARE yields the equivalent CARE
			\begin{equation}
				\tilde A(\mu)^\top Y+Y\tilde A(\mu)-Y\tilde G(\mu)Y+Q(\mu)=0,
			\end{equation}
			where $\tilde A=E^{-1}A$, $\tilde G=E^{-1}BB^\top E^{-\top}$, $Q=C^\top C$. By (A2) and Lemma~\ref{lem:inv_lipschitz}, the maps $\tilde A,\tilde G,Q$ are Lipschitz on $\mathcal P$, with constants $L_A,L_G,L_Q$ in Frobenius norm. The original solution is recovered by
			\begin{equation}
				X^*(\mu)=E(\mu)^{-\top}Y(\mu)E(\mu)^{-1}.
			\end{equation}
			
			\paragraph{Step 2: Continuity of $Y(\mu)$} 
			Fix an arbitrary $\mu_0\in\mathcal P$. Let $Y_0:=Y(\mu_0)$ and define $\Phi_0=\tilde A(\mu_0)-\tilde G(\mu_0)Y_0$. Define the Lyapunov operator $\mathcal L_0(W)=\Phi_0^\top W+W\Phi_0$ and constants
			\begin{align}
				l_0 &= \|\mathcal L_0^{-1}\|^{-1},\\
				p_0 &= \sup_{\|W\|_F=1}\|\mathcal L_0^{-1}(Y_0W+W^\top Y_0)\|_F,\\
				q_0 &= \sup_{\|W\|_F=1}\|\mathcal L_0^{-1}(Y_0WY_0)\|_F.
			\end{align}
			For perturbations $\Delta A=\tilde A(\mu)-\tilde A(\mu_0)$, $\Delta G=\tilde G(\mu)-\tilde G(\mu_0)$, $\Delta Q=Q(\mu)-Q(\mu_0)$, set
			\[
			\delta=\|\Delta A\|_F+\|\Delta G\|_F\|Y_0\|_2,\
			\hat g=\|\tilde G(\mu_0)\|_2+\|\Delta G\|_2,
			\]
			and
			\[
			\epsilon=\frac1{l_0}\|\Delta Q\|_F+p_0\|\Delta A\|_F+q_0\|\Delta G\|_F.
			\]
			Since the coefficient matrices $\tilde A$, $\tilde G$, and $Q$ are Lipschitz continuous, we have
			\[
			\|\Delta A\|_F \le L_A\|\mu-\mu_0\|_2,\
			\|\Delta G\|_F \le L_G\|\mu-\mu_0\|_2,\
			\|\Delta Q\|_F \le L_Q\|\mu-\mu_0\|_2.
			\]
			Consequently, there exists $\eta_{\mu_0}>0$ such that whenever $\|\mu-\mu_0\|_2<\eta_{\mu_0}$,
			\begin{equation}
				\delta+\sqrt{l_0\,\hat g\,\epsilon}<\frac{l_0}{2}. \label{eq:first_hyp}
			\end{equation}
			This verifies the hypothesis of Sun's theorem. Applying it yields a local Lipschitz estimate
			\begin{equation}
				\|Y(\mu)-Y(\mu_0)\|_F \le C_{Y,\mu_0}\|\mu-\mu_0\|_2. \label{eq:first_lip}
			\end{equation}
			Since $\mu_0$ was arbitrary, $Y$ is continuous on $\mathcal P$.
			
			\paragraph{Step 3: Uniform bounds}
            As $Y$ is continuous on compact $\mathcal{P}$, there exists $M_Y > 0$ with
           \begin{equation}
              \|Y(\mu)\|_F \le M_Y, \quad \forall \mu \in \mathcal{P}. \label{eq:A8}
           \end{equation}

           The closed-loop matrix $\Phi(\mu) = \tilde{A}(\mu) - \tilde{G}(\mu)Y(\mu)$ is stable for each $\mu \in \mathcal{P}$. Since the mapping $\mu \mapsto \Phi(\mu)$ is continuous, the associated Lyapunov operator 
           \begin{equation*}
               \mathcal{L}_\mu(W) = \Phi(\mu)^\top W + W \Phi(\mu)
           \end{equation*}
           is also continuous with respect to $\mu$. Because $\Phi(\mu)$ is stable, $\mathcal{L}_\mu$ is invertible, which implies $\sigma_{\min}(\mathcal{L}_\mu) > 0$ for every $\mu \in \mathcal{P}$. 

           Since the minimum singular value $\sigma_{\min}(\cdot)$ is a continuous function of the matrix entries, the mapping $\mu \mapsto \sigma_{\min}(\mathcal{L}_\mu)$ is continuous. As $\mathcal{P}$ is compact, this continuous function attains its minimum over $\mathcal{P}$, yielding a uniform positive bound
           \begin{equation}
               l_{\min} := \min_{\mu \in \mathcal{P}} \sigma_{\min}(\mathcal{L}_\mu) > 0. \label{eq:A9}
           \end{equation}

           Consequently, the operator norm of the inverse satisfies
           \begin{equation*}
               \|\mathcal{L}_\mu^{-1}\| = \frac{1}{\sigma_{\min}(\mathcal{L}_\mu)} \le \frac{1}{l_{\min}}, \quad \forall \mu \in \mathcal{P}.
           \end{equation*}
           Equivalently, $\|\mathcal{L}_\mu^{-1}\|^{-1} \ge l_{\min}$.

           Consequently, the quantities
           \begin{align}
                p(\mu) &:= \sup_{\|W\|_F = 1} \|\mathcal{L}_\mu^{-1}(Y(\mu)W + W^\top Y(\mu))\|_F, \label{eq:p_mu} \\
               q(\mu) &:= \sup_{\|W\|_F = 1} \|\mathcal{L}_\mu^{-1}(Y(\mu)W Y(\mu))\|_F, \label{eq:q_mu}
            \end{align}
            satisfy the uniform bounds
           \begin{equation}
               p(\mu) \le p_{\max} := \frac{2 M_Y}{l_{\min}}, \quad q(\mu) \le q_{\max} := \frac{M_Y^2}{l_{\min}}. \label{eq:A10}
           \end{equation}
			
			\paragraph{Step 4: Local Lipschitz continuity of $Y(\mu)$} 
			Fix $\mu_0\in\mathcal P$ and use the uniform constants $l_{\min},p_{\max},q_{\max},M_Y$ to define
			\[
			\delta=\|\Delta A\|_F+\|\Delta G\|_F M_Y,\
			\epsilon=\frac1{l_{\min}}\|\Delta Q\|_F+p_{\max}\|\Delta A\|_F+q_{\max}\|\Delta G\|_F.
			\]
			By Lipschitzness, there exist $C_\delta,C_\epsilon>0$ independent of $\mu_0$ with $\delta\le C_\delta\|\mu-\mu_0\|_2$ and $\epsilon\le C_\epsilon\|\mu-\mu_0\|_2$. Also $\hat g$ is bounded by $g_{\max}+L_G\|\mu-\mu_0\|_2$, where $g_{\max}=\max_\mu\|\tilde G(\mu)\|_2$. Hence there exists $\eta>0$ (independent of $\mu_0$) such that whenever $\|\mu-\mu_0\|_2<\eta$,
			\begin{equation}
				\delta+\sqrt{l_{\min}\,\hat g\,\epsilon}<\frac{l_{\min}}{2}. \label{eq:second_hyp}
			\end{equation}
			Applying Sun's theorem again gives a constant $C_Y>0$, independent of $\mu_0$, such that
			\begin{equation}
				\|Y(\mu)-Y(\mu_0)\|_F \le C_Y\|\mu-\mu_0\|_2,\ \|\mu-\mu_0\|_2<\eta. \label{eq:second_lip}
			\end{equation}
			
			\paragraph{Step 5: Global Lipschitz continuity of $Y(\mu)$} 
			For arbitrary $\mu,\nu\in\mathcal P$, convexity ensures the segment $[\mu,\nu]\subset\mathcal P$. Choose $m$ such that $\|\mu-\nu\|_2/m<\eta$, and partition the segment into points $\mu_j=\mu+\frac{j}{m}(\nu-\mu)$. Applying the local estimate \eqref{eq:second_lip} on each subinterval and summing via the triangle inequality gives
			\begin{equation}
				\|Y(\mu)-Y(\nu)\|_F \le C_Y\|\mu-\nu\|_2, \label{eq:Y_global}
			\end{equation}
			so $Y$ is globally Lipschitz.
			
			\paragraph{Step 6: Global Lipschitz continuity of $X^*(\mu)$}
			\begin{align*}
				X^*(\mu)-X^*(\nu)
				&= E(\mu)^{-\top}Y(\mu)E(\mu)^{-1} - E(\nu)^{-\top}Y(\nu)E(\nu)^{-1} \\
				&= E(\mu)^{-\top}\bigl(Y(\mu)-Y(\nu)\bigr)E(\mu)^{-1} \\
				&\quad + E(\mu)^{-\top}Y(\nu)\bigl(E(\mu)^{-1}-E(\nu)^{-1}\bigr) \\
				&\quad + \bigl(E(\mu)^{-\top}-E(\nu)^{-\top}\bigr)Y(\nu)E(\nu)^{-1}.
			\end{align*}
			Taking the Frobenius norm and using its submultiplicativity, we obtain
			\begin{align*}
				\|X^*(\mu)-X^*(\nu)\|_F
				&\le \|E(\mu)^{-\top}\|_2\,\|Y(\mu)-Y(\nu)\|_F\,\|E(\mu)^{-1}\|_2 \\
				&\quad + \|E(\mu)^{-\top}\|_2\,\|Y(\nu)\|_F\,\|E(\mu)^{-1}-E(\nu)^{-1}\|_F \\
				&\quad + \|E(\mu)^{-\top}-E(\nu)^{-\top}\|_F\,\|Y(\nu)\|_F\,\|E(\nu)^{-1}\|_2.
			\end{align*}
			Applying the uniform bounds and Lipschitz estimates, we can obtain
			\[
			\|X^*(\mu)-X^*(\nu)\|_F
			\le \left(M_E^2 C_Y + 2M_E^3 L_E M_Y\right)\|\mu-\nu\|_2.
			\]
			Thus $X^*$ is globally Lipschitz on $\mathcal P$. This completes the proof.
	\end{proof}
		
	\bibliographystyle{siamplain}
	\bibliography{references}
\end{document}